\documentclass[aop,MSNbibl,seceqn,dvips]{arximspdf}

\doi{10.1214/13-AOP834} %kopijuoti is PTS
\volume{41}
\issue{6}
\pubyear{2013}
\firstpage{3973}
\lastpage{4001}

\makeatletter
\newcommand{\onrarrow}[1]{\mathrel{\vbox{\m@th\ialign{##\crcr
  $\hfil\scriptstyle\ #1\ \ \hfil$\crcr\noalign{\kern0.5pt\nointerlineskip}%
  \rightarrowfill\crcr}}}}

\newcommand{\rrvert}{\vert}
\newcommand{\llvert}{\vert}
\newcommand{\eqref}[1]{(\ref{#1})}
\newtheorem{thmm}{Theorem}[section]
\newtheorem{cor}[thmm]{Corollary}

\newtheorem{lem}[thmm]{Lemma}
\newtheorem{prop}[thmm]{Proposition}

\newproclaim{defn}[thmm]{Definition}

\newproclaim{exep}[thmm]{Example}

\newproclaim{rem}[thmm]{Remark}

\newcommand{\R}{\mathbb{R}}
\newcommand{\N}{\mathbb{N}}

\newcommand{\essinf}{\operatorname{ess inf}}

\makeatother

\begin{document}
\begin{frontmatter}

\title{Minimal supersolutions of convex BSDEs}
\runtitle{Minimal supersolutions of convex BSDEs}

\begin{aug}
\author{\fnms{Samuel} \snm{Drapeau}\thanksref{t2}\ead[label=e1]{drapeau@math.hu-berlin.de}},
\author{\fnms{Gregor} \snm{Heyne}\corref{}\thanksref{t1}\ead[label=e2]{heyne@math.hu-berlin.de}}
\and
\author{\fnms{Michael} \snm{Kupper}\thanksref{t2}\ead[label=e3]{kupper@math.hu-berlin.de}}
\thankstext{t2}{Supported from MATHEON project E.11.}
\thankstext{t1}{Supported from MATHEON project E.2.}
\runauthor{S. Drapeau, G. Heyne and M. Kupper}
\affiliation{Humboldt University Berlin}
\address{Humboldt University Berlin\\
Unter den Linden 6\\
D-10099 Berlin\\
Germany\\
\printead{e1}\\
\phantom{E-mail:\ }\printead*{e2}\\
\phantom{E-mail:\ }\printead*{e3}} %adresu isvedimo komanda gale!
\end{aug}

% HISTORY:
\received{\smonth{8} \syear{2011}}
\revised{\smonth{1} \syear{2013}}

% ABSTRACT
%
\begin{abstract}
We study the nonlinear operator of mapping the terminal value $\xi$ to
the corresponding minimal supersolution of a backward stochastic
differential equation with the generator being monotone in $y$, convex in
$z$, jointly lower semicontinuous and bounded below by an affine function
of the control variable $z$.
We show existence, uniqueness, monotone convergence, Fatou's lemma and
lower semicontinuity of this operator.
We provide a comparison principle for minimal supersolutions of BSDEs.
\end{abstract}

% KEYWORDS
\end{frontmatter}

%s1 #&#
\section{Introduction}

On a filtered probability space, where the filtration is generated by a
$d$-dimensional Brownian motion $W$, we consider the process $\hat
{\mathcal{E}}^{g}(\xi)$ given by
\[
\hat{\mathcal{E}}^{g}_{t}(\xi)=\essinf\bigl\{
Y_t \in L^0_t\dvtx ( Y,Z )\in\mathcal{A}(\xi,g)
\bigr\},\qquad t \in [ 0,T ],
\]
where $\mathcal{A}(\xi,g)$ is the set of all pairs of c\`adl\`ag
\emph
{value processes} $Y$ and \emph{control processes} $Z$ such that
%
%e1.1 #&#
%
\begin{equation}
\label{eqintrodyn} Y_s-\int_s^t
g_{u}(Y_u,Z_u)\,du+\int_s^t
Z_{u}\,dW_{u}\geq Y_t \quad\mbox {and}\quad Y_T
\geq\xi
\end{equation}
for all $0\leq s\leq t\leq T$.
Here the \emph{terminal condition} $\xi$ is a random variable, the
\emph
{generator} $g$ a measurable function of $(y,z)$ and the pair $(Y,Z)$
is a \emph{supersolution} of the \emph{backward stochastic differential
equation \textup{(}BSDE\textup{)}},
\[
dY_t=-g_t(Y_t,Z_t)\,dt+Z_t\,dW_t,\qquad
t\in[0,T); Y_T=\xi.
\]
%
% \eqref{eqintrodyn}.

The main objective of this paper is to state conditions which guarantee
that there exists a unique minimal supersolution.
More precisely, we show that the process $\mathcal{E}^{g}(\xi)=\lim_{s\downarrow\cdot, s\in\mathbb{Q}}\hat{\mathcal{E}}^{g}_s(\xi)$
is a
modification of $\hat{\mathcal{E}}^{g}(\xi)$ and equals the value
process of the unique minimal supersolution, that is, there exists a
unique control process $\hat Z$ such that $(\mathcal{E}^{g}(\xi),\hat
Z)\in\mathcal{A}(\xi,g)$.
The existence theorem immediately yields a comparison theorem for
minimal supersolutions.
We also study the stability of the minimal supersolution with respect
to the terminal condition and the generator.
We show that the mapping $\xi\mapsto\hat{\mathcal{E}}^{g}_{0}(\xi)$
is a nonlinear expectation, fulfills a monotone convergence theorem and
Fatou's lemma on the same domain as the expectation operator $E[\cdot
]$, and consequently is ${L^1}$-lower semicontinuous.

%Before we give further details on the specific choice of our spaces
%and assumptions, let us recall that related problems have been
%investigated throughout the literature before.
%Most notably, nonlinear expectations have been a prominent topic in
%mathematical economics since Allais famous paradox, see \cite[Section~2.2]{foellmer01}.
Nonlinear expectations have been a prominent topic in mathematical
economics since Allais's famous paradox, see F\"{o}llmer and Schied
\cite{foellmer01}, Section~2.2.
Typical examples are the monetary risk measures introduced by Artzner
et~al. \cite{artzner01} and F\"{o}llmer and Schied \cite{foellmer02}, Peng's $g$
and $G$-expectations
\cite{peng03,peng07,peng08}, the variational preferences by
Maccheroni, Marinacci and
Rustichini \cite
{marinacci01} and the recursive utilities by Duffie and Epstein \cite
{epstein03}.
Especially the $g$-expectation, which is defined as the initial value
of the solution of a BSDE, is closely related to $\mathcal
{E}^{g}_{0}(\cdot)$, since each pair $(Y,Z)$ that solves the BSDE
corresponding to \eqref{eqintrodyn} is also a supersolution and hence
an element of $\mathcal{A}(\xi,g)$.
The concept of a supersolution of a BSDE appears already in El~Karoui,
Peng and Quenez \cite{karoui01}, Section~2.2.
For further references see Peng \cite{peng99}, who derives monotonic limit
theorems for supersolutions of BSDEs and proves the existence of a
minimal constrained supersolution.
%Another related concept are stochastic target problems, which were
%introduced and studied by \cite{ST2002}, by means of controlled
%stochastic differential equations and dynamic programming methods.

Our first contribution is to provide a setting where we relax the usual
Lipschitz requirements for the generator $g$.
Namely, we suppose that $g$ is convex with respect to $z$, monotone in
$y$, jointly lower semicontinuous, and bounded below by an affine
function of the control variable $z$.
To see in an intuitive way the role these assumptions play in deriving
the existence and uniqueness of a control process $\hat Z$ such that
$(\mathcal{E}^{g}(\xi),\hat Z)\in\mathcal{A}(\xi,g)$, let us suppose
for the moment that $g$ is positive.
%In order to provide an intuition on how these assumptions contribute
%toward the existence and uniqueness of a control process $\hat Z$ such
%that $(\mathcal{E}^{g}(\xi),\hat Z)\in\mathcal{A}(\xi,g)$, let us
%suppose for the moment that $g$ is positive.
Given an adequately good space of control processes, the value process
of each supersolution and the process $\hat{\mathcal{E}}^g(\xi)$ are in
fact supermartingales.
By suitable pasting, we may now construct a decreasing sequence $(Y^n)$
of supersolutions, whose pointwise limit is again a supermartingale and
equal to $\hat{\mathcal{E}}^g(\xi)$ on all dyadic rationals.
%The candidate value process is then given by the c\`adl\`ag
%supermartingale $\mathcal{E}^{g}(\xi)$, the standard right limit of $
Since the generator $g$ is positive, it can be shown that $\mathcal
{E}^{g}(\xi)$ lies below $\hat{\mathcal{E}}^g(\xi)$, $P$-almost surely,
at any time.
This suggests to consider the c\`adl\`ag supermartingale $\mathcal
{E}^{g}(\xi)$ as a candidate for the value process of the minimal
supersolution.
However, it is not clear a priori that the sequence $(Y^n)$ converges
to $\mathcal{E}^{g}(\xi)$ in some suitable sense.
Yet, taking into account the additional supermartingale structure, in
particular the Doob--Meyer decomposition, it follows that $(Y^n)$
converges $P\otimes dt$-almost surely to $\mathcal{E}^{g}(\xi)$.
It remains to obtain a unique control process $\hat Z$ such that
$(\mathcal{E}^{g}(\xi),\hat Z)\in\mathcal{A}(\xi,g)$.
To that end, we prove that, for monotone sequences of supersolutions, a
positive generator yields, after suitable stopping, a uniform
$L^1$-bound for the sequence of supremum processes of the associated
sequence of stochastic integrals.
This, along with a result by Delbaen and Schachermayer \cite{DS94},
and standard compactness
arguments and diagonalization techniques yield the candidate control
process $\hat Z$ as the limit of a sequence of convex combinations.
Now, joint lower semicontinuity of $g$, positivity and convexity in $z$
allow us to use Fatou's lemma to verify that the candidate processes
$(\mathcal{E}^{g}(\xi),\hat Z)$ are a supersolution of the BSDE.
Thus, $\mathcal{E}^{g}(\xi)$ is in fact the value process of the
minimal supersolution and a modification of $\hat{\mathcal
{E}}^{g}(\xi)$.
Finally, the uniqueness of $\hat Z$ follows from the uniqueness of the
Doob--Meyer decomposition of the c\`adl\`ag supermartingale~$\mathcal
{E}^{g}(\xi)$.

Let us give further reference of related assumptions and methods in the
existing literature.
Delbaen, Hu and Bao \cite{BDH10} consider superquadratic BSDEs with
generators that are
positive and convex in $z$ but do not depend on~$y$.
However, their principal aim and method differ from ours.
Indeed, they primarily study the well-posedness of superquadratic BSDEs
by establishing a dual link between cash additive time-consistent
dynamic utility functions and supersolutions of BSDEs.
To view supersolutions as supermartingales is one of the key ideas in
our approach, and we make ample use of the rich structure
supermartingales provide.
The classical limit theorem of supermartingales has been used by
El~Karoui and Quenez \cite
{EQ95} in the theory of BSDEs, when studying the problem of option
pricing in incomplete financial markets.
However, the analysis is done via dual formulations and only for linear
generators that do not depend on $y$.
The construction of solutions of BSDEs by monotone approximations is
also a classical tool; see, for example, Kobylanski \cite
{kobylanski01} for
quadratic generators and Briand and Hu \cite{briand02} for generators
that are in
addition convex in $z$.
The application of compactness theorems such as Delbaen and
Schachermayer \cite{DS94}, Lemma
A1.1 or \cite{delbaen03}, Theorem A, in order to derive
existence of BSDEs seems to be new to the best of our knowledge.
Often existence proofs rely on a priori estimates combined with a fixed
point theorem (see, e.g., \cite{karoui01}) or on constructing Cauchy
sequences in complete spaces; see, for example, Briand and Confortola
\cite{briand01} or
Ankirchner, Imkeller and
Dos~Reis \cite{imkeller01}.
Recent exceptions are R\'eveillac \cite{Reveillac11} and Heyne, Kupper
and Mainberger \cite{HKM1011} who use a
compactness result given in Barlow and Protter \cite{BarProt}.
As already mentioned, Peng \cite{peng99} studies the existence and
uniqueness of minimal supersolutions.
However, he assumes a Lipschitz generator, a square integrable terminal
condition, and employs a very different approach.
It is based on a monotonic limit theorem, \cite{peng99}, Theorem~2.4
and the penalization method introduced in El~Karoui et~al. \cite
{EKKPPQ97}, and it
leads to increasing sequences of supersolutions.
Parallel to us, Cheridito and Stadje \cite{cs11} have investigated
existence and stability
of supersolutions of BSDEs.
They consider generators that are convex in $z$ and Lipschitz in $y$.
However, their setting and methods are quite different from ours.
Namely, they approximate by discrete time BSDEs and work with terminal
conditions that are bounded lower semicontinuous functions of the
Brownian motion.
An interesting equivalence between the minimal supersolution and the
solution of a reflected BSDEs is given in Peng and Xu \cite{PengXu}.
In \cite{HKM1011} the authors show the existence of the minimal
supersolution for generators that are lower semicontinuous, monotone in
the value variable, are bounded below by an affine function of the
control variable and which satisfy a specific normalization condition.
Finally, given our local $L^1$-bounds, the compactness underlying the
construction of the candidate control process is a special case of
results obtained by Delbaen and Schachermayer \cite{delbaen03}.

Our second contribution is to allow for \emph{local supersolutions},
that is, for supersolutions $(Y,Z)$, where the stochastic integral of
$Z$ is only a local martingale.
However, in order to avoid so-called ``doubling strategies,'' present
even for the simplest generator $g\equiv0$ (see Dudley \cite
{dudley02} or
Harrison and Pliska \cite{harrison01}, Section~6.1), we require in
addition that $\int
Z\,dW$ is a supermartingale.
This specification interacts nicely with a positive generator and
happens to be particularly adequate in establishing stability
properties of the minimal supersolution with respect to the terminal
condition or the generator.
In particular, it allows us to formulate theorems such as monotone
convergence and Fatou's lemma for the nonlinear operator $\hat
{\mathcal
{E}}^g_{0}(\cdot)$ on the same domain as the standard expectation
$E[\cdot]$ and to obtain its $L^1$-lower semicontinuity.
Moreover, under some additional integrability on the terminal
condition, our approach also allows us to derive existence results with
control processes, whose stochastic integrals belong to $\mathcal{H}^1$.

Dropping the positivity assumption, the value and control processes of
our supersolutions are supermartingales under another measure closely
linked to the generator $g$.
In fact, for a positive generator we have supermartingales with respect
to the initial probability measure $P$, while for a nonpositive
generator, which is bounded below by an affine function of the control
variable, we consider supermartingales under the measure given by the
corresponding Girsanov transform.
%This yields a generator dependent concept of admissibility.
%The implication thereof is illustrated by giving a minimal
%supersolution based approach to the well known problem of exponential
%expected utility maximization, where this admissibility is related in
%a natural way to the market price of risk.

The paper is organized as follows.
In Section~\ref{secsetting} we fix our notation and the setting.
We define minimal supersolutions and introduce our main conditions and
structural properties of $\hat{\mathcal{E}}^g(\xi)$ in Section~\ref{sectarpointro}.
Finally in Section~\ref{secL1target} we state and prove our main
results, existence and stability theorems.
%which concludes with an example on maximizing expected exponential
%utility.

%s2 #&#
\section{Setting and notation}\label{secsetting}
We consider a fixed time horizon $T>0$ and a filtered probability space
$(\Omega,\mathcal{F}, (\mathcal{F}_t )_{t \in [ 0,T
]},P)$, where the filtration $ ( \mathcal{F}_t )$ is generated
by a $d$-dimensional Brownian motion $W$ and fulfills the usual conditions.
We further assume that $\mathcal{F}=\mathcal{F}_{T}$.
The set of $\mathcal{F}$-measurable and $\mathcal{F}_t$-measurable
random variables is denoted by $L^0$ and $L^0_t$, respectively, where
random variables are identified in the $P$-almost sure sense.
The sets $L^p$ and $L^p_t$ denote the set of random variables in $L^0$
and $L^0_t$, respectively, with finite $p$-norm, for $p \in [
1,+\infty ]$.
Throughout this work, inequalities and strict inequalities between any
two random variables or processes $X^1,X^2$ are understood in the
$P$-almost sure or in the $P\otimes dt$-almost sure sense,
respectively; that is, $X^1\leq(<) X^2$ is equivalent to $P [
X^1\leq(<) X^2  ]=1$ or $P\otimes dt  [ X^1\leq(<) X^2
]=1$, respectively.
Given a process $X$ and $t\in[0,T]$ we denote $X^\ast_{t}:=\sup_{s
\in
[ 0,t  ]}\vert X_s\vert$.
We denote by $\mathcal{T}$ the set of stopping times with values in
$[0,T]$ and hereby call an increasing sequence of stopping times $(\tau
^n)$, such that $P[\bigcup_{n}\{ \tau^n=T \} ]=1$, a \emph{localizing
sequence of stopping times}.
By $\mathcal{S}:=\mathcal{S} ( \R )$ we denote the set of all
c\`adl\`ag progressively measurable processes $Y$ with values in $\R$.
For $p \in [ 1,+\infty [$, we further denote by $\mathcal
{L}^p:=\mathcal{L}^p ( W  )$ the set of progressively
measurable processes $Z$ with values in $\R^{1\times d}$, such that
$\Vert Z\Vert_{\mathcal{L}^p}:=E[(\int_{0}^{T}Z_sZ_s^{\top}
\,ds)^{p/2}]^{1/p}< +\infty$.
For any $Z \in\mathcal{L}^p$, the stochastic integral $(\int_0^t Z_s
\,dW_s)_{t \in [ 0,T  ]}$ is well defined (see \cite
{protter01}) and is by means of the Burkholder--Davis--Gundy
inequality, a continuous martingale.
For the $\mathcal{L}^p$-norm, the set $\mathcal{L}^p$ is a Banach
space; see \cite{protter01}.
We further denote by $\mathcal{L}:=\mathcal{L} ( W  )$ the set
of progressively measurable processes with values in $\R^{1\times d}$,
such that there exists a localising sequence of stopping times $(\tau
^n)$ with $Z1_{ [0,\tau^n  ]} \in\mathcal{L}^1$, for all
$n\in
\N$.
Here again, the stochastic integral $\int_{}^{}Z\,dW$ is well defined and
is a continuous local martingale.

For adequate integrands $a,Z$, we generically write $\int_{}^{}a \,ds$ or
$\int_{}^{}Z \,dW$ for the respective integral processes $(\int_{0}^{t}a_s \,ds)_{t \in [ 0,T  ]}$ and $( \int_{0}^{t}Z_s \,dW_s
)_{t \in [ 0,T  ]}$.
Finally, given a sequence $(x_n)$ in some convex set, we say that a
sequence $(y_n)$ is in the \emph{asymptotic convex hull} of $(x_n)$, if
$y_n \in \operatorname{conv}\{ x_{n},x_{n+1},\ldots \} $, for all $n$.

A \emph{generator} is a jointly measurable function g from $\Omega
\times[0, T ]\times\mathbb{R}\times\mathbb{R}^{1\times d}$ to
$\mathbb{R}\cup\{+\infty\}$ where $\Omega\times[0, T ]$ is endowed
with the progressive $\sigma$-field.

%
%%
%
%s3 #&#
\section{Minimal supersolutions of BSDEs}\label{sectarpointro}

%s3.1 #&#
\subsection{Definitions}

Given a generator $g$, and a \emph{terminal condition} $\xi\in L^0$, a
pair $(Y,Z)\in\mathcal{S}\times\mathcal{L}$ is a \emph{supersolution}
of the BSDE
\[
dY_t=-g_t(Y_t,Z_t)\,dt+Z_t\,dW_t,\qquad
t\in[0,T); Y_T=\xi,
\]
if, for all $s,t\in[0,T]$, with $s\leq t$, it holds
%
%e3.1 #&#
%
\begin{equation}
\label{eqcentralineq} Y_s-\int_s^t
g_{u}(Y_u,Z_u)\,du+\int_s^t
Z_{u}\,dW_{u}\geq Y_t\quad \mbox {and}\quad Y_T
\geq\xi.
\end{equation}
For such a supersolution $(Y,Z)$, we call $Y$ the \emph{value process}
and $Z$ its \emph{control process}.
Due to the c\`adl\`ag property, relation \eqref{eqcentralineq} holds
for all stopping times $0\leq\sigma\leq\tau\leq T$, in place of $s$
and $t$, respectively.
Note that the formulation in \eqref{eqcentralineq} is equivalent to
the existence of a c\`adl\`ag increasing process $K$, with $K_{0}=0$,
such that
%
%e3.2 #&#
%
\begin{equation}
\label{eqstandardnotationsupersolution}\qquad Y_{t}=\xi+ \int
_t^T g_{u}(Y_u,Z_u)\,du+(K_{T}-K_{t})-
\int_t^T Z_{u}\,dW_{u},\qquad t
\in[0,T].
\end{equation}
Although the notation in \eqref{eqstandardnotationsupersolution} is
standard in the literature concerning supersolutions of BSDEs (see,
e.g., \cite{karoui01,peng99}), we will keep with \eqref
{eqcentralineq} since the proofs of our main results exploit this structure.
We consider only those supersolutions $(Y,Z)\in\mathcal{S}\times
\mathcal{L}$ of a BSDE where $Z$ is \emph{admissible}, that is, where
the continuous local martingale $\int Z\,dW$ is a supermartingale.
We are then interested in the set
%
%e3.3 #&#
%
\begin{equation}
\label{setA} \mathcal{A}(\xi,g)= \bigl\{ (Y,Z) \in\mathcal{S}\times\mathcal{L} \dvtx Z \mbox{ is admissible and } \eqref{eqcentralineq} \mbox{ holds} \bigr\}
\end{equation}
and the process
%
%e3.4 #&#
%
\begin{equation}
\label{eqdefiphi} \hat{\mathcal{E}}^{g}_{t}(\xi)=\essinf
\bigl\{ Y_t \in L^0_t\dvtx ( Y,Z )\in\mathcal{A}(
\xi,g) \bigr\},\qquad t \in [ 0,T ].
\end{equation}
By $\hat{\mathcal{E}}^g$ we mean the functional mapping terminal
conditions $\xi\in L^0$ to the process $\hat{\mathcal{E}}^g (
\xi
)$.
%if there is no ambiguity about $\xi$.
If necessary, we write $\mathcal{A}_T (\xi,g  )$ and $\hat
{\mathcal{E}}^g_{\cdot,T} ( \xi )$ for $\mathcal{A}
( \xi
,g  )$ and $\hat{\mathcal{E}}^g ( \xi )$, respectively,
to indicate their dependence on the time horizon.
%Since the set $\mathcal{A}\left( \xi,g \right)$ and therefore $\hat{
%indicate this by writing $\mathcal{A}_T\left( \xi,g \right)$ and $\hat{
Note that the essential infima in \eqref{eqdefiphi} can be taken over
those $(Y,Z)\in\mathcal{A}(\xi,g)$, where $Y_{T}=\xi$.
A pair $(Y,Z)$ is called a \emph{minimal supersolution}, if $(Y,Z)\in
\mathcal{A}(\xi,g)$, and if for any other supersolution $(Y',Z')\in
\mathcal{A}(\xi,g)$, holds $Y_{t}\leq Y'_{t}$, for all $t\in[0,T]$.

%
%s3.2 #&#
\subsection{\texorpdfstring{General properties of $\mathcal{A}(\cdot,g)$ and $\hat{\mathcal{E}}^g$}
{General properties of A(.,g) and E g}}\label{secpropertiesAE}

In this section we collect various statements regarding the properties
of $\mathcal{A} ( \cdot,g  )$ and $\hat{\mathcal{E}}^g$.
The first lemma ensures that the set of admissible control processes is
stable under pasting and that we may concatenate elements of $\mathcal
{A}(\xi,g)$ along stopping times and partitions of our probability space.
%
%le3.1 #&#
\begin{lem}\label{lemadmpasting}
Fix a generator $g$, a terminal condition $\xi\in L^0$, a stopping
time $\sigma\in\mathcal{T}$ and $(B^n) \subset\mathcal{F}_\sigma$ a
partition of $\Omega$:
\begin{longlist}[(1)]
\item[(1)] Let $(Z^n)\subset\mathcal L$ be admissible.
Then $\bar Z = Z^11_{[0,\sigma]} + \sum_{n\ge1} Z^n1_{B_n}1_{]\sigma
,T]}$ is admissible.
\item[(2)] Let $((Y^n,Z^n))\subset\mathcal{A} (\xi,g)$ such that
$Y^1_{\sigma}1_{B^n}\geq Y^n_{\sigma}1_{B^n}$, for all $n\in\N$.
Then $(\bar Y,\bar Z)\in\mathcal{A} ( \xi,g  )$, where
%
%e3.5 #&#
%
\begin{eqnarray}
\label{lemadmpastingpartition} \bar Y&=&Y^11_{[0,\sigma[}+\sum
_{n\ge1} Y^n 1_{B^n}
1_{[\sigma,T]}\quad \mbox{and}
\nonumber
\\[-8pt]
\\[-8pt]
\nonumber
 \bar Z&=&Z^11_{[0,\sigma]}+\sum
_{n\ge1} Z^n 1_{B^n} 1_{]\sigma,T]}.
\end{eqnarray}
\end{longlist}
\end{lem}
\begin{pf}
(1) Let $M^n$ and $\bar M$ denote the stochastic integrals of the
$Z^n$ and $\bar Z$, respectively.
It follows from $(Z^n)\subset\mathcal L$ and from $(B_n)$ being a
partition that $\bar Z\in\mathcal L$ and that $\int_{s\vee\sigma
}^{t\vee
\sigma}\bar Z_u\,dW_u=\sum1_{B_n}\int_{s\vee\sigma}^{t\vee\sigma
}Z^n_u\,dW_u$.\vadjust{\goodbreak}
Now observe that the admissibility of all $Z^n$ yields
\begin{eqnarray*}
E [\bar M_t-\bar M_s| \mathcal{F}_s ]&=&E
\bigl[M^1_{(t\wedge
\sigma)\vee s}-M^1_s|
\mathcal{F}_{s} \bigr]
\\
&&{}+ E \biggl[\sum_{n\ge
1}1_{B_n}E
\bigl[M^n_{t\vee\sigma}-M^n_{s\vee\sigma}| \mathcal
{F}_{s\vee\sigma} \bigr]| \mathcal{F}_s \biggr]\le0
\end{eqnarray*}
for $0\le s\le t\le T$.

(2) $\bar{Z}$ is admissible by item (1).
Since $Y^1_{\sigma}1_{B^n}\geq Y^n_{\sigma}1_{B^n}$, for all $n\in\N$,
it follows on the set $\{ s<\sigma\leq t \} $ that
\begin{eqnarray*}
&&Y_s^1-\int_s^\sigma
g_u\bigl(Y^1_u,Z^1_u
\bigr)\,du+\int_s^\sigma Z^1_u
\,dW_u-\int_\sigma^t g_u(
\bar Y_u,\bar Z_u)\,du+\int_\sigma^t
\bar Z_u\,dW_u
\\
&&\qquad\geq Y_\sigma^1-\sum_{n\geq1}1_{B^n}
\biggl(\int_\sigma^t g_u
\bigl(Y^n_u,Z^n_u\bigr)\,du-\int
_\sigma^t Z^n_u\,dW_u
\biggr)
\\
&&\qquad\geq\sum_{n\geq1}1_{B^n}
\biggl(Y_{\sigma}^n-\int_\sigma^t
g_u\bigl(Y^n_u,Z^n_u
\bigr)\,du+\int_\sigma^t Z^n_u\,dW_u
\biggr)\geq\sum_{n\geq
1}1_{B^n}Y^n_t.
\end{eqnarray*}
Hence,
\begin{eqnarray*}
&&\bar{Y}_s-\int_s^t
g_u(\bar{Y}_u,\bar{Z}_u)\,du+\int
_s^t\bar{Z}_u\,dW_u
\\
&&\qquad\geq1_{\{ \sigma> t \} } Y_{t}^1+\sum
_{n\geq1}1_{B^n} \bigl( 1_{\{ \sigma\leq s \} }
Y_{t}^n+ 1_{\{ s<\sigma\leq t \}
}Y_{t}^n
\bigr)=\bar Y_{t}
\end{eqnarray*}
and thus $(\bar Y, \bar Z)\in\mathcal{A}(\xi,g)$.
\end{pf}
For convenience, a generator is said to be:
\begin{longlist}[\textsc{(Dec)}]%[label=,leftmargin=40pt]
\item[\textsc{(Pos)}] positive if $g ( y,z  )\geq0$, for all
$(y,z)\in\R\times\R^{1\times d}$.
\item[\textsc{(Inc)}] increasing if $g (y,z  )\geq g (
y^\prime,z )$, for all $y,y'\in\R$ with $y\geq y^\prime$, and all
$z \in\R^{1\times d}$.
\item[\textsc{(Dec)}] decreasing if $g (y,z  )\leq g (
y^\prime,z )$, for all $y,y'\in\R$ with $y\geq y^\prime$, and all
$z \in\R^{1\times d}$.
\end{longlist}

%The following Lemma states that under the assumption of a positive
%generator the value process of a supersolution is a supermartingale.
%To view supersolutions as supermartingales is one of the key ideas in
%our approach.
In the following lemma, we show that the value process of a
supersolution is a supermartingale if the generator is positive.
%
%le3.2 #&#
\begin{lem}\label{remdecompY}
Let $g$ be a generator fulfilling \textsc{(Pos)}, and $\xi\in L^0$ be a
terminal condition such that $\xi^- \in L^1$.
Let $(Y,Z)\in\mathcal{A}(\xi,g)$.
Then $\xi\in L^1$, $Y$ is a supermartingale, $Z$ is unique and $Y$ has
the unique decomposition
%
%e3.6 #&#
%
\begin{equation}
\label{eqdecompY} Y=Y_{0}-A+M,
\end{equation}
where $M$ denotes the supermartingale $\int Z\,dW$, and $A$ is a
predictable, increasing, c\`adl\`ag process with $A_{0}=0$.
\end{lem}
\begin{pf}
Relation \eqref{eqcentralineq}, positivity of $g$, admissibility of
$Z$ and $\xi^- \in L^1$ imply $E[\vert Y_t\vert]<+\infty$, for all
$t \in[0,T]$.
Since $-\xi^-\leq\xi\leq Y_T $, we deduce that $\xi\in L^1$.
Again, from \eqref{eqcentralineq}, admissibility of $Z$ and
positivity of $g$ we derive by taking conditional expectation, that
$Y_{s}\geq E[Y_{t}| \mathcal{F}_{s}]$.
Thus $Y$ is a supermartingale with $Y_t\geq E [ \xi| \mathcal
{F}_t  ]$.
Relation \eqref{eqcentralineq} implies further that there exist an
increasing and c\`adl\`ag process $K$, with $K_{0}=0$, such that \eqref
{eqdecompY} holds with $A=\int g(Y,Z)\,ds+K$.
Note that $A$ is optional and therefore predictable due to the Brownian
filtration; see \cite{Revuz1999}, Corollary V.3.3.
Since $Y$ is a c\`adl\`ag supermartingale the Doob--Meyer theorem, see
\cite{protter01}, Theorem III.3.13, implies the unique decomposition
$Y=Y_{0}+\tilde M-\tilde A$, where $\tilde M$ is a local martingale and
$\tilde A$ is an increasing process which is predictable, and $\tilde
M_{0}=\tilde A_{0}=0$.
In our filtration every local martingale is continuous (see \cite
{protter01}, Corollary IV.3.1, page 187), and thus $\tilde A$ is c\`
adl\`ag.
Hence $A$ and $\tilde A$ and $M$ and $\tilde M$ are indistinguishable.
Moreover, from the predictable representation property of local
martingales and from $P(\bigcup_{n}\{\tau_{n}=T\})=1$, for $\tau
^n=\inf
\{t\geq0| |M_{t}|\geq n \}\wedge T$, we obtain the $P\otimes
dt$-almost sure uniqueness of $Z$.
\end{pf}
%
%The next proposition addresses the dependence of $\mathcal{A}\left(
%The first two properties are crucial in the proof of the existence and
%uniqueness theorem in Section~\ref{secL1target}.
%The third item concerns the monotonicity of $\hat{\mathcal{E}}^g\left(
%Combined with the existence theorem, this yields in fact a comparison
%principle for minimal supersolutions of BSDEs.
%Finally, the last item concerns the cash (super/sub) additivity of the
%functional $\hat{\mathcal{E}}^g(\xi)$.
%
%pr3.3 #&#
\begin{prop}\label{propproperties1}
For $t \in [ 0,T  ]$, generators $g,g^\prime$ and terminal
conditions $\xi,\xi^\prime\in L^0$, it holds:
\begin{longlist}[(1)]
\item[(1)] the set $\{ Y_{t}\dvtx (Y,Z)\in\mathcal{A}(\xi,g) \} $ is directed
downwards;
\item[(2)] assuming \textsc{(Pos)}, $\xi^-\in L^1$
and $\mathcal{A} ( \xi,g  )\neq\varnothing$, then for all
$\varepsilon>0$, there exists $ ( Y^\varepsilon,Z^\varepsilon
) \in\mathcal{A} ( \xi,g  )$ such that $\hat{\mathcal
{E}}_t^g ( \xi )\geq Y^\varepsilon_t-\varepsilon$;
\item[(3)] \emph{(monotonicity)} if
$\xi
^\prime\leq\xi$ and $g^\prime ( y,z  )\leq g ( y,z
)$, for all $y,z \in\R\times\R^{1\times d}$, then $\mathcal
{A} ( \xi^\prime,g^\prime )\supset\mathcal{A} (
\xi,g
)$ and $\hat{\mathcal{E}}_t^{g^\prime} ( \xi^\prime
)\leq\hat{\mathcal{E}}_t^g ( \xi )$;
\item[(4)] \emph{(convexity)} if $(y,z)\mapsto g(y,z)$ is jointly
convex, then $\mathcal{A} ( \lambda\xi+ ( 1-\lambda
)\xi
^\prime,g )\supset\lambda\mathcal{A} ( \xi,g
)+ (
1-\lambda )\mathcal{A} ( \xi^\prime,g  )$, for all
$\lambda\in ( 0,1  )$, and so
\[
\hat{\mathcal{E}}_{t}^g \bigl( \lambda\xi+ ( 1-\lambda )
\xi ^\prime \bigr)\leq\lambda\hat{\mathcal{E}}^g_t
( \xi )+ ( 1-\lambda )\hat{\mathcal{E}}^g_t \bigl(
\xi^\prime \bigr).
\]
\item[(5)] for $m \in L^{0}_t$:
\begin{itemize}
\item\emph{(cash superadditivity)} assuming \textsc{(Inc)} and $m\geq
0$, then $\hat{\mathcal{E}}_t^g  ( \xi+m  )\geq\hat
{\mathcal
{E}}_t^g ( \xi )+m$.
\item\emph{(cash subadditivity)} assuming \textsc{(Dec)}, $m\geq0$, and
the existence of $(Y,Z)\in\mathcal{A}(\xi,g)$, such that $\mathcal
{A}_{t}(Y_{t}+m,g)\neq\varnothing$, then $\hat{\mathcal
{E}}^g_{t}
(\xi+m )\leq\hat{\mathcal{E}}^g_{t} (\xi )+m$.
\item%\label{abcd} \label{itemtransinv}
\emph{(cash additivity)}
assuming that $g$ does not depend on $y$, the existence of $(Y,Z)\in
\mathcal{A}(\xi,g)$, such that $\mathcal{A}_{t}(Y_{t}+m^+,g)\neq
\varnothing$, and the existence of $(Y,Z)\in\mathcal{A}(\xi+m,g)$, such
that $\mathcal{A}_{t}(Y_{t}+m^-,g)\neq\varnothing$, then $\hat
{\mathcal
{E}}^g_{t} (\xi+m )=\hat{\mathcal{E}}^g_{t} (\xi
)+m$.
\end{itemize}
\end{longlist}
\end{prop}
\begin{pf}
(1) Given $(Y^{i},Z^{i})\in\mathcal{A}(\xi,g)$, for $i=1,2$, we have
to construct $(\bar Y,\bar Z)\in\mathcal{A}(\xi,g)$, such that $\bar
Y_{t}\leq\min\{ Y^{1}_{t},Y^{2}_{t} \} $.
To this end, we define the stopping time
\[
\tau=\inf\bigl\{ s> t\dvtx Y^{1}_{s}> Y^{2}_{s}
\bigr\} \wedge T
\]
and set\vspace*{1pt} $\bar Y= Y^{1} 1_{[0,\tau[}+Y^{2} 1_{[\tau,T[}$, $ \bar
Y_T=\xi
$, and $\bar Z = Z^{1} 1_{[0,\tau]}+Z^{2} 1_{]\tau,T]}$.
Since $Y_{\tau}^1\geq Y^2_{\tau}$, Lemma~\ref{lemadmpasting} yields
$(\bar Y,\bar Z)\in\mathcal{A}(\xi,g)$ and by definition holds $\bar
Y_{t}=\min\{ Y^{1}_{t},Y^{2}_{t} \} $.

(2) In view of the first assertion, there exists a sequence
$((\tilde Y^n,\tilde Z^n))\subset\mathcal{A}(\xi,g)$ such that
$(\tilde
Y^n_t)$ decreases to $\hat{\mathcal{E}}^g_t (\xi )$.
Set $Y^n=\tilde Y^11_{[0,t)}+\tilde Y^n1_{[t,T]}$ and $Z^n=\tilde
Z^11_{[0,t]}+\tilde Z^n1_{(t,T]}$.
From Lemma~\ref{lemadmpasting} it follows that $((Y^n,Z^n))\subset
\mathcal{A}(\xi,g)$ and $(Y^n_t)$ decreases to $\hat{\mathcal
{E}}^g_t (\xi )$ by construction.
% Relation \eqref{eqcentralineq}, $Y_0^n = Y_{0}^1$, and $g$ positive
%yield
% \[
% \int_0^t Z^n_s dW_s\ge\xi-\int_t^T Z^n_s dW_s+\int_{0}^{T}g_s\left(
%Y_s^n,Z_s^n \right)\,ds-Y_0^n\geq-\xi^--\int_t^T Z^n_s dW_s-Y_{0}^1.
% \]
% Taking conditional expectation with respect to $\mathcal{F}_t$ and
%using the supermartingal property of $\int_{}^{}Z^n dW$ yield
% \[
% \int_0^t Z^n_s dW_s\ge-E\left[\xi^-\Mid\mathcal{F}_t\right]-Y_{0}^1.
% \]
Lemma~\ref{remdecompY} implies that $\hat{\mathcal{E}}^g_t (\xi
)\geq E[\xi^-|\mathcal{F}_t]$.
Hence, given $\varepsilon>0$, the sets $B^n=A^n\setminus A^{n-1}\in
\mathcal{F}_t$, where $A^n=\{\hat{\mathcal{E}}^g_t (\xi
)\ge
Y^n_t-\varepsilon\}$ and $A^0=\varnothing$, form a partition of $\Omega$.
Since $(Y^n_t)$ is decreasing, it follows that $Y^1_t1_{B^n}\geq
Y^n_t1_{B^n}$, for all $n\in\N$.
Consequently, by means of Lemma~\ref{lemadmpasting}, $(\bar Y, \bar
Z)$, defined as in \eqref{lemadmpastingpartition}, is an element of
$\mathcal{A}(\xi,g)$ and $\hat{\mathcal{E}}^g_t (\xi
)\ge\bar
Y_t-\varepsilon$ by construction.

(3) Follows from definitions \eqref{setA} and \eqref{eqdefiphi}.
%%shows that $(Y,Z)\in\mathcal{A}\left( \xi,g \right)$, and relation
%

(4) The joint convexity of $g$ yields $(\lambda Y +(1-\lambda
)Y',\lambda Z +(1-\lambda)Z')\in\mathcal{A}(\lambda\xi+(1-\lambda
)\xi
',g)$, for all $(Y,Z)\in\mathcal{A}(\xi,g)$, all $(Y',Z')\in
\mathcal
{A}(\xi',g)$ and all $\lambda\in(0,1)$.
Hence, $\lambda\mathcal{A} ( \xi,g  )+ ( 1-\lambda
)\mathcal{A} ( \xi^\prime,g  )\subset\mathcal{A} (
\lambda
\xi+ ( 1-\lambda )\xi^\prime,g  )$ and in particular,
$\hat{\mathcal{E}}^g_{t}(\lambda\xi+(1-\lambda)\xi^\prime)\leq
\lambda
\hat{\mathcal{E}}^g_{t}(\xi)+(1-\lambda)\hat{\mathcal
{E}}^g_{t}(\xi
^\prime)$.

(5) Let us show the cash superadditivity.
For $m \in L^{0}_t$ with $m\geq0$, given $(Y,Z)\in\mathcal{A}(\xi
+m,g)$, and $0\leq s\leq t^\prime\leq T$, it follows from \eqref
{eqcentralineq} and \textsc{(Inc)} that
\begin{eqnarray*}
&&Y_{s}-m1_{ [ t,T  ]} ( s )-\int_{s}^{t^\prime
}g_{u}
\bigl(Y_{u}-m1_{[t,T]}(u),Z_{u}\bigr)\,du+\int
_{s}^{t^\prime}Z_{u}\,dW_{u}
\\
&&\qquad\geq Y_{s}-m1_{ [ t,T  ]} ( s )-\int_{s}^{t^\prime
}g_{u}(Y_{u},Z_{u})\,du+
\int_{s}^{t^\prime}Z_{u}\,dW_{u}\geq
Y_{t^\prime
}-m1_{ [ t,T  ]} \bigl( t^\prime \bigr).
\end{eqnarray*}
Hence, $(Y-m1_{[t,T]},Z)\in\mathcal{A}(\xi,g)$ and thus $\hat
{\mathcal
{E}}^g_{t}(\xi+m)-m \geq\hat{\mathcal{E}}^g_{t}(\xi)$.
For the cash subadditivity the same argument yields
\[
Y_{s}+m1_{ [ t,T  ]} ( s )-\int_{s}^{t^\prime
}g_{u}
\bigl(Y_{u}+m1_{[t,T]}(u),Z_{u}\bigr)\,du+\int
_{s}^{t^\prime
}Z_{u}\,dW_{u}\geq
Y_{t^\prime}+m1_{ [ t,T  ]} \bigl( t^\prime \bigr)
\]
for all $t\leq s\leq t^\prime\leq T$, and all $(Y,Z)\in\mathcal
{A}(\xi,g)$.
In order to apply our usual pasting argument we now need the
assumption that $\mathcal{A}_{t}(Y_{t}+m,g)\neq\varnothing$.\vspace*{1pt}
It provides $(\tilde Y,\tilde Z)\in\mathcal{A}_{t}(Y_{t}+m,g)$ such
that we may construct $(\bar Y, \bar Z)\in\mathcal{A}(\xi+m,g)$, with
$Y_{t}+m=\bar Y_{t}$ and thus $\hat{\mathcal{E}}^g_{t}(\xi)+m \geq
\hat
{\mathcal{E}}^g_{t}(\xi+m)$.
The cash additivity in the case where $g$ is independent of $y$
follows from $\hat{\mathcal{E}}^g_{t}(\xi)+m=\hat{\mathcal
{E}}^g_{t}(\xi
+m^+)-m^-=\hat{\mathcal{E}}^g_{t}(\xi+m+m^-)-m^-=\hat{\mathcal
{E}}^g_{t}(\xi+m)$, since \textsc{(Dec)} and \textsc{(Inc)} are
simultaneously fulfilled.
\end{pf}
Proposition~\ref{propproperties1} addresses the dependence of $\mathcal{A} ( \xi,g
)$ on $\xi$ and $g$ and its impact on $\hat{\mathcal
{E}}^g (
\xi )$. The first two assertions are crucial in the subsequent
proof of the existence and uniqueness theorem in Section~\ref{secL1target}.
The third assertion concerns the monotonicity of $\hat
{\mathcal{E}}^g ( \xi )$ with respect\vadjust{\goodbreak} to $\xi$ and $g$.
Combined with the existence theorem, this yields, in fact, a comparison
principle for minimal supersolutions of BSDEs. The last assertion
concerns the cash (super/sub) additivity of the functional $\hat
{\mathcal{E}}^g ( \xi )$.

We now prove that for a positive generator $\hat{\mathcal{E}}^g
(\xi
)$ is in fact a supermartingale, which, in addition, dominates
its right-hand-limit process.
This is crucial for the proof of the existence and uniqueness theorem.
%
%pr3.4 #&#
\begin{prop}\label{propsinnlosphisubmartingale}
Let $g$ be a generator fulfilling \textsc{(Pos)}, and $\xi\in L^0$ be a
terminal condition such that $\xi^- \in L^1$.
Suppose that $\mathcal{A} ( \xi,g  )\neq\varnothing$. Then
$\hat{\mathcal{E}}^g (\xi )$ is a supermartingale,
\[
\mathcal{E}^{g}_s(\xi):=\lim_{t\downarrow s,t\in\mathbb{Q}}
\hat {\mathcal {E}}^g_t(\xi) \qquad\mbox{for all } s\in[0,T),\qquad
\mathcal {E}^{g}_T(\xi):=\xi,
\]
is a well-defined c\`adl\`ag supermartingale and
%
%e3.7 #&#
%
\begin{equation}
\label{eqphiincreasing} \hat{\mathcal{E}}^g_s (\xi )
\geq\mathcal{E}^{g}_s(\xi) \qquad\mbox{for all } s\in[0,T].
\end{equation}
\end{prop}
\begin{pf}
Note first that $\hat{\mathcal{E}}^g(\xi)$ is adapted by definition.
Furthermore, given $(Y,Z) \in\mathcal{A}(\xi,g)\neq\varnothing$, Lemma~\ref{remdecompY} implies $\xi\in L^1$ and $Y_t\geq E [ \xi|
\mathcal{F}_t  ]$.
Hence $Y_{t}\geq\hat{\mathcal{E}}^g_{t}(\xi)\geq E[\xi| \mathcal
{F}_{t}]$ and $\hat{\mathcal{E}}^g_{t}(\xi)\in L^1$.
As for the supermartingale property and~\eqref{eqphiincreasing}, fix
$0\leq s\leq t\leq T$.
In view of item (2) of Proposition~\ref
{propproperties1}, for all $\varepsilon>0$, there exists
$(Y^\varepsilon,Z^\varepsilon)\in\mathcal{A}(\xi,g)$ such that
$\hat{\mathcal{E}}^g_s(\xi)\ge Y^\varepsilon_s-\varepsilon$.
Due to \eqref{eqcentralineq} it follows
%
%e3.8 #&#
%e3.9 #&#
%
\begin{eqnarray}
\label{schwachsinn3} \hat{\mathcal{E}}^g_t(\xi) &\le&
Y^\varepsilon_t\le Y^\varepsilon _s-\int
_s^t g_{u}\bigl(Y^\varepsilon_u,Z^\varepsilon_u
\bigr) \,du + \int_s^t Z^\varepsilon
_u \,dW_u
\nonumber\\
&\le&\hat{\mathcal{E}}^g_s ( \xi )-\int
_s^t g_{u}\bigl(Y^\varepsilon_u,Z^\varepsilon_u
\bigr) \,du + \int_s^t Z^\varepsilon_u
\,dW_u+\varepsilon\\
&\le&\hat{\mathcal{E}}^g_s(
\xi)+ \int_s^tZ^\varepsilon _u
\,dW_u+\varepsilon.\nonumber
\end{eqnarray}
Taking conditional expectation on both sides with respect to $\mathcal
{F}_s$, the supermartingale property of $\int Z^\varepsilon \,dW$ yields
$\hat{\mathcal{E}}^g_s(\xi)\ge E[\hat{\mathcal{E}}_t^g(\xi)|
\mathcal
{F}_s]$, and so $\hat{\mathcal{E}}^g ( \xi )$ is a
supermartingale.
That $\mathcal{E}^{g}(\xi)$ is a well-defined c\`adl\`ag
supermartingale follows from Karatzas and Shreve \cite{karatzas01},
Proposition~1.3.14.
Finally, \eqref{eqphiincreasing} follows directly from~\eqref
{schwachsinn3} and the definition of $\mathcal{E}^g ( \xi )$.
\end{pf}
%
%re3.5 #&#
\begin{rem}\label{remEhutundmodication}
The previous proposition suggests to consider the c\`adl\`ag
supermartingale $\mathcal{E}^{g}(\xi)$ as a candidate for the value
process of the minimal supersolution.
Note further that if $\mathcal{E}^g ( \xi )$ is the value
process of the minimal supersolution, it is a modification of $\hat
{\mathcal{E}}^g ( \xi )$ by definition.
%Hence, in the following, whenever we assume that $\mathcal{E}^g\left(
%regard to fixed deterministic times, we may write either of them.
\end{rem}
The final result of this section shows that our setup allows us to
derive various properties that are important in the context of
nonlinear expectations and dynamic risk measures.
In particular,\vadjust{\goodbreak} we prove that $\mathcal{E}^g ( \xi )$, if it
is the value process of the minimal supersolution, fulfills the
flow-property and, under the additional assumption $g(y,0)=0$, for all
$y\in\R$, we show projectivity, with time-consistency as a special case.
In the context of BSDE solutions such properties were first established
in \cite{peng03}, for the case of Lipschitz generators.
For dynamic risk measures the (strong) time-consistency has been
investigated in discrete time in \cite{cheridito01,foellmer03} as
well as in continuous time in \cite{bion01,delbaen05}, for instance.
%
%pr3.6 #&#
\begin{prop}\label{propsinnlos}
For $t \in [ 0,T  ]$, generator $g$ and terminal condition
$\xi\in L^0$, it holds:
\begin{longlist}[(1)]
\item[(1)] $\hat{\mathcal{E}}^g_{s,T} (
\xi )\leq\hat{\mathcal{E}}^g_{s,t}( \hat{\mathcal
{E}}^g_{t,T} ( \xi ) )$, for all $0\leq s \leq t$.
Suppose that $\mathcal{E}^g ( \xi )$ is a minimal
supersolution, then the \emph{flow-property} holds; that is,
%
%e3.10 #&#
%
\begin{equation}
\mathcal{E}^g_{s,T} ( \xi )= \mathcal{E}^g_{s,t}
\bigl( \mathcal {E}^g_{t,T} ( \xi )\bigr) \qquad\mbox{for all }0
\leq s\leq t.
\end{equation}

\item[(2)] If $g ( y,0  )=0$, for
all $y \in\R$, then $\hat{\mathcal{E}}^g_{s} ( \hat{\mathcal
{E}}^g_{t} ( \xi )  )\leq\hat{\mathcal
{E}}^g_{s} (
\xi )$, for all $0\leq s\leq t$.
Assuming \textsc{(Pos)}, $\xi^-\in L^1$, and supposing that $\mathcal
{E}^g ( \xi )$ is a minimal supersolution, then $\mathcal
{E}^g ( \xi )$ is \emph{time-consistent}, that is,
%
%e3.11 #&#
%
\begin{equation}
\mathcal{E}^g_{s} \bigl( \mathcal{E}^g_{t}
( \xi ) \bigr)= \mathcal{E}^g_{s} ( \xi )\qquad \mbox{for all
}0\leq s\leq t.
\end{equation}
\item[(3)] Assuming \textsc{(Pos)},
$g
( y,0  )=0$, for all $y \in\R$, $\xi^-\in L^1$, and $\mathcal
{E}^g ( \xi )$ is a minimal supersolution, then the \emph
{projectivity} holds, that is,
%
%e3.12 #&#
%
\begin{equation}
\label{eqmoinblablabla} \mathcal{E}^g_{s} \bigl(
1_A\mathcal{E}^g_{t} ( \xi ) \bigr)=
\mathcal{E}^g_{s} ( 1_A\xi ) \qquad\mbox{for all
}0\leq s\leq t \mbox{ and } A \in\mathcal{F}_t.
\end{equation}
\end{longlist}
\end{prop}
\begin{pf}
(1) Fix $0\leq s\leq t$.
Obviously, $ ( Y_s,Z_s  )_{s \in [ 0,t  ] }\in
\mathcal{A}_{t}( \hat{\mathcal{E}}^g_{t,T} ( \xi ),g)$, for
all $  ( Y,Z  ) \in\mathcal{A}_T ( \xi,g  )$.
Hence $\hat{\mathcal{E}}_{s,t}^g( \hat{\mathcal{E}}_{t,T}^g (
\xi
))\leq\hat{\mathcal{E}}_{s,T}^g ( \xi )$.
Suppose now that $\mathcal{E}_{\cdot,T}^g ( \xi )$ is a
minimal supersolution with corresponding admissible control process
$\hat{Z}\in\mathcal{L}$.
For all $(Y,Z) \in\mathcal{A}_{t}( \mathcal{E}_{t,T}^g ( \xi
),g)$, holds $Y_t\geq\mathcal{E}_{t,T}^g ( \xi )$ and,
with the same argumentation as in Lemma~\ref{lemadmpasting}, we can
paste in a monotone way to show that $(\bar{Y},\bar{Z})\in\mathcal
{A}_T ( \xi,g  )$, where $\bar{Y}=Y1_{ [ 0,t
[}+\mathcal{E}_{\cdot,T}^g ( \xi )1_{ [ t,T
]}$ and
$\bar{Z}=Z1_{ [ 0,t  ]}+\hat{Z}1_{ ]t,T  ]}$.
Thus, by definition, $\mathcal{E}^g_{s,t}( \mathcal{E}^g_{t,T} (
\xi ))\geq\mathcal{E}^g_{s,T} ( \xi )$.

(2)
Given $ ( Y,Z  )\in\mathcal{A} ( \xi,g  )$, we
define $\bar{Y}=Y1_{ [ 0,t  [}+\hat{\mathcal
{E}}^g_{t} (
\xi )1_{ [ t,T  ]}$ and $\bar{Z}=Z1_{ [ 0,t
]}$.
Since $Y_t\geq\hat{\mathcal{E}}_{t}^g ( \xi )$ and
$g
(y,0  )=0$, it is straightforward to verify that $(\bar{Y},\bar
{Z})\in\mathcal{A}( \hat{\mathcal{E}}_t^g ( \xi ),g)$.
From $Y_s\geq\bar{Y}_s$, for all $s \in [ 0,t  ]$, follows
that $\hat{\mathcal{E}}_{s}^g( \hat{\mathcal{E}}_{t}^g ( \xi
))\leq\hat{\mathcal{E}}_{s}^g ( \xi )$, for all $s \in
[ 0,t  ]$.
The case where $\mathcal{E}^g ( \xi )$ is a minimal
supersolution and assumption~\textsc{(Pos)} holds, follows from \eqref
{eqmoinblablabla} for $A=\Omega$.

(3)
Fix $A \in\mathcal{F}_t$.
Suppose that $\mathcal{E}^g ( \xi )$ is a minimal
supersolution\vspace*{1pt} with corresponding control process $\hat{Z}$.
Then, from $\xi^-\in L^1$ and Lemma~\ref{remdecompY} it follows that
$\mathcal{E}^g ( \xi )$ is a supermartingale and $\xi\in L^1$.

Given $ ( Y,Z  )\in\mathcal{A} ( 1_{A}\mathcal
{E}_t^g(\xi
),g  )$, it follows from $(1_{A}\mathcal{E}_t^g(\xi))^-\in L^1$
and Lem\-ma~\ref{remdecompY} that $Y_t\geq E [ 1_{A}\mathcal
{E}^g_t ( \xi )| \mathcal{F}_t  ]=1_A\mathcal
{E}_t^g ( \xi )$.
Since $g ( y,0  )=0$, it is straightforward to check that
$\tilde{Y}=Y1_{ [ 0,t  [}+\mathcal{E}^g_t ( \xi
)
1_A 1_{ [ t,T  ]}$, and $\tilde{Z}=Z1_{ [ 0,t
]}$ is
such that $( \tilde{Y},\tilde{Z})\in\mathcal{A} (
1_{A}\mathcal
{E}_{t}^g ( \xi ),g  )$.
We can henceforth assume that $Y_{s}=1_{A}\mathcal{E}_{t}^g ( \xi
)$, for all $s \geq t$.
Now, we define $\bar{Y}=Y 1_{[0,t[}+\mathcal{E}^g ( \xi
)1_{A}1_{ [ t,T  ]} $ and $\bar{Z}=Z1_{[0,t]}+\hat
{Z}1_{A}1_{]t,T]} $,
for $ 0\leq s<t\leq t^\prime\leq T$ holds
\begin{eqnarray*}
&&\bar{Y}_{s} -\int_{s}^{t^\prime}g (
\bar{Y}_{u},\bar{Z}_{u} )\,du+\int_{s}^{t^\prime}
\bar{Z}_u\,dW_u
\\
&&\qquad\geq Y_{t}+ \biggl(-\int_{t}^{t^\prime}g_{u}
\bigl(\mathcal {E}^g_{u} ( \xi ),\hat{Z}_u
\bigr)\,du+\int_{t}^{t^\prime}\hat {Z}_u\,dW_u
\biggr)1_A
\\
&&\qquad\geq \biggl(\mathcal{E}^g_{t} ( \xi )-\int
_{t}^{t^\prime
}g_{u} \bigl(
\mathcal{E}^g_{u} ( \xi ),\hat{Z}_u \bigr)\,du+
\int_{t}^{t^\prime}\hat{Z}_u\,dW_u
\biggr)1_A\geq\mathcal {E}^g_{t^\prime} ( \xi
)1_A.
\end{eqnarray*}
Hence, for all $0\leq s\leq t^\prime\leq T$ it holds that
\[
\bar{Y}_{s} -\int_{s}^{t^\prime}g (
\bar{Y}_{u},\bar{Z}_{u} )\,du+\int_{s}^{t^\prime}
\bar{Z}_u \,dW_u\geq Y_{t^\prime}1_{\{
t^\prime < t \} }+
\mathcal{E}^g_{t^\prime} ( \xi ) 1_{A}1_{\{ t\leq t^\prime \} }=
\bar{Y}_{t^\prime}
\]
and $\bar{Y}_{T}=1_{A}\xi$, which implies that $(\bar{Y},\bar
{Z})\in
\mathcal{A} ( 1_{A}\xi,g  )$.
Since $\bar{Y}_s=Y_s$, for all $s \leq t$, we deduce $\mathcal
{E}^g_{s} ( 1_A \xi )\leq\mathcal{E}^g_{s} (
1_{A}\mathcal{E}_{t} ( \xi )  )$.

On the other hand, consider $ ( Y,Z  )\in\mathcal{A} (
1_A \xi,g  )$.
From $Y_{t}\geq E [1_{A}\xi| \mathcal{F}_t  ]=1_{A}E
[\xi| \mathcal{F}_t  ]$, we obtain $Y_{t}1_{A^c}\geq0$.
Since $\mathcal{E}^g ( \xi )$ is a minimal supersolution, it
follows that $Y_t\geq\mathcal{E}^g_t ( \xi )1_{A}$.
Indeed, let $B=\{ Y_t<\mathcal{E}^g_{t} ( \xi )1_A \} $; then
$Y_{t}1_{A^c}\geq0$ implies $B\subset A$.
Consequently, by arguments similar to those in Lemma~\ref
{lemadmpasting}, the processes $\tilde{Y}=\mathcal{E}^g ( \xi
)(1_{[0,t[}+1_{B^c}1_{[t,T]})+Y1_B1_{[t,T]}$ and $\tilde
{Z}=\hat
{Z}(1_{[0,t[}+1_{B^c}1_{[t,T]})+Z1_B1_{[t,T]}$ are such that $(\tilde
{Y},\tilde{Z})\in\mathcal{A} ( \xi,g  )$, which implies
$P [ B  ]=0$.
It is also straightforward to check that $\tilde{Y}=Y1_{ [ 0,t
[}+\mathcal{E}^g ( \xi ) 1_A 1_{ [ t,T  ]}$
and $\tilde{Z}=Z1_{ [ 0,t  ]}+\hat Z1_{ ( t,T
]}1_{A}$ are such that $( \tilde{Y},\tilde{Z})\in\mathcal{A} (
1_{A}\xi,g  )$.
Thus we can assume that $Y_{t}=1_{A}\mathcal{E}_{t}^g ( \xi
)$.
Defining $\bar{Y}=Y 1_{[0,t[}+\mathcal{E}^g_t ( \xi
)1_{A}1_{ [ t,T  ]} $ and $\bar{Z}=Z1_{[0,t]}$, it holds $(
\bar{Y},\bar{Z} )\in\mathcal{A} ( 1_A\mathcal{E}_t^g (
\xi
),g  )$.
Thus $\mathcal{E}^g_s ( 1_{A}\mathcal{E}^g_t ( \xi )
)\leq\mathcal{E}^g_s ( 1_A\xi )$, since $\bar
{Y}_s=Y_s$, for all $s \leq t$.
\end{pf}

%
%s4 #&#
\section{Existence, uniqueness and stability}\label{secL1target}
In this section, we give conditions that guarantee the existence and
uniqueness of a minimal supersolution.
We show that the corresponding value process is given by $\mathcal
{E}^g(\xi)$.
Moreover, we analyze the stability of $\hat{\mathcal{E}}^g(\xi)$ with
respect to perturbations of the terminal condition or the generator.
In addition to the assumptions \textsc{(Pos)} and \textsc{(Inc)} or \textsc{(Dec)} introduced above, we require convexity of $g$ in the control
variable and joint lower semicontinuity.
To that end, we say that a generator $g$ is:
\begin{longlist}%[label=,leftmargin=40pt]
\item[\textsc{(Con)}] convex if $g ( y,\lambda z+ ( 1-\lambda
)z^\prime )\leq\lambda g ( y,z  )+ (
1-\lambda )g ( y,z^\prime )$, for all $y \in\R$, all
$z,z^\prime\in\R^{1 \times d}$ and all $\lambda\in ( 0,1
)$;
\item[\textsc{(Lsc)}] if $(y,z)\mapsto g(y,z)$ is lower semicontinuous.
\end{longlist}

%s4.1 #&#
\subsection{Existence and uniqueness of minimal supersolutions}\label{secex}
The following theorem on existence and uniqueness of a minimal
supersolution is the first main result of this paper.
%
%th4.1 #&#
\begin{thmm}\label{thmmexistence}
Let $g$ be a generator fulfilling \textsc{(Pos)}, \textsc{(Lsc)}, \textsc{(Con)} and either \textsc{(Inc)} or \textsc{(Dec)} and $\xi\in L^0$ be a
terminal condition, such that $\xi^-\in L^1$.
If $\mathcal{A}(\xi,g) \neq\varnothing$, then there exists a unique
minimal supersolution $(\hat Y,\hat Z)\in\mathcal{A}(\xi,g)$.
Moreover, $\mathcal{E}^g ( \xi )$ is the value process of
the minimal supersolution, that is, $(\mathcal{E}^g ( \xi
),\hat{Z})\in\mathcal{A}(\xi,g)$.
\end{thmm}
%
% \textcolor{red}{A study of the proof shows that we do not need the
%generator to be jointly lower semicontinuous, if we assume
% However, for the study of the stability of $\mathcal{E}^g\left( \xi
%the generator, see section \ref{secstab}, we will need that property.
%}
% %\end{rem}
Note that, under the assumptions of Theorem~\ref{thmmexistence},
Remark~\ref{remEhutundmodication} implies that the process $\mathcal
{E}^g ( \xi )$ is a modification of $\hat{\mathcal
{E}}^g
( \xi )$.
Further, in the context of finding minimal elements in some set, the
assumption $\mathcal{A}(\xi,g) \neq\varnothing$ is quite standard; see
\cite{peng99} for an example in the setting of minimal supersolutions.
However, let us point out that in many applications $\mathcal{A}(\xi,g)
\neq\varnothing$ might be guaranteed by specific model assumptions; see,
for instance, an example on utility maximization in Heyne \cite{heyne12}.
It might also be automatically granted under further assumptions (see
Cheridito and Stadje \cite{cs11}) or, for instance, if the BSDE
$Y_t-\int_{t}^{T}g_s
(Y_s,Z_s  )\,ds+\int_{t}^{T}Z_s \,dW_s=\hat\xi$ has a solution
$ (
Y,Z  )\in\mathcal{S}\times\mathcal{L}$, such that $Z$ is admissible.
In the latter case, $\mathcal{A} ( \xi,g  )\neq\varnothing$,
for all $\xi\in L^0$ such that $\xi^- \in L^1$, with $\hat\xi\ge
\xi$.

\begin{pf*}{Proof of Theorem \ref{thmmexistence}}
\textit{Step} 1. \textit{Uniqueness.}
Given $\hat Z\in\mathcal{L}$ such that $(\mathcal{E}^g(\xi),\hat
Z)\in
\mathcal{A}(\xi,g)$, the definition of $\mathcal{E}^g(\xi)$ implies
that for any other supersolution $(Y, Z')\in\mathcal{A}(\xi,g)$ holds
$\mathcal{E}^g_{t}(\xi)\leq Y_{t}$, for all $t\in[0,T]$.
The uniqueness of $\hat Z$ follows as in Lemma~\ref{remdecompY}.

The remainder of the proof provides existence of $\hat Z\in\mathcal{L}$
such that $(\mathcal{E}^g(\xi),\break \hat Z)\in\mathcal{A}(\xi,g)$.

\textit{Step} 2. \textit{Construction of an approximating sequence.}
For any $n,i\in{\mathbb N}$, let $t^n_i=iT/2^n$.
There exist $((Y^n,Z^n))\subset\mathcal{A}(\xi,g)$ such that
%
%e4.1 #&#
%
\begin{equation}
\label{eqstep1nprocess} \hat{\mathcal{E}}^g_{t^n_i}(\xi)\ge
Y^n_{t^n_i}-1/n\qquad \mbox{for all } n\in\N \mbox{ and all }
i=0,\ldots, 2^n-1
\end{equation}
and
%
%e4.2 #&#
%
\begin{equation}
\label{eqthmmexynmono} Y^{n}_{t}\geq Y^{n+1}_{t}\qquad
\mbox{for all } t\in[0,T]\mbox{ and all } n\in\mathbb{N}.
\end{equation}
Indeed, by means of Proposition~\ref{propproperties1}(2), for each $n\in\N$, we may select a family
$((Y^{n,i},Z^{n,i}))_{i=0,\ldots, 2^n-1}$ in $ \mathcal{A}(\xi,g)$, such
that $\hat{\mathcal{E}}^g_{t^n_i}(\xi)\ge Y^{n,i}_{t^n_i}-1/n$,
$i=0,\ldots, 2^n-1$.
We suitably paste this family in order to obtain \eqref{eqstep1nprocess}.
We start with
\[
\bar Y^{n,0} = Y^{n,0}, \qquad\bar Z^{n,0}=Z^{n,0}
\]
and continue by recursively setting, for $i=1,\ldots, 2^n-1$,
\begin{eqnarray*}
\bar Y^{n,i} &= &\bar Y^{n,i-1}1_{[0,\tau^n_{i}[} +
Y^{n,i} 1_{[\tau
^n_{i},T[},\qquad \bar Y^{n,i}_{T}=\xi,
\\
\bar Z^{n,i} &=& \bar Z^{n,i-1}1_{[0,\tau^n_{i}]} +
Z^{n,i} 1_{]\tau^n_{i},T]},
\end{eqnarray*}
where $\tau^n_{i}$ are stopping times given by $\tau^n_{i}=\inf\{t>
t^n_i \dvtx \bar Y^{n,i-1}_{t}> Y^{n,i}_{t}\}\wedge T$.
From the definition of these stopping times and Lemma~\ref
{lemadmpasting} follows that the pairs $(\bar Y^{n,i},\bar Z^{n,i})$,
$i=0,\ldots, 2^n-1$, are elements of $\mathcal{A}(\xi,g)$.
Hence the sequence
\[
\bigl(\bigl( Y^n,Z^n\bigr):=\bigl(\bar Y^{n,2^n-1},
\bar Z^{n,2^n-1}\bigr)\bigr)
\]
fulfills \eqref{eqstep1nprocess} by construction.
Note that $(( Y^n,Z^n))$ is not necessarily monotone in the sense of
\eqref{eqthmmexynmono}.
However, this can be achieved by pasting similarly.
More precisely, we choose
\[
\bar Y^{1} = Y^{1}, \qquad \bar Z^{1}=Z^{1},
\]
and continue by recursively setting, for $n\in\mathbb{N}$,
\begin{eqnarray*}
\bar Y^n &=& \sum_{i=0}^{2^n-1}
\bigl(Y^n 1_{[t^{n}_i,\tau
^{n}_i[}+\bar Y^{n-1} 1_{[\tau^{n}_i,t^{n}_{i+1}[}
\bigr), \qquad\bar Y_T^n=\xi,
\\
\bar Z^n &= &\sum_{i=0}^{2^n-1}
\bigl(Z^n 1_{]t^{n}_i,\tau
^{n}_i]}+\bar {Z}^{n-1} 1_{]\tau^{n}_i,t^{n}_{i+1}]}
\bigr),
\end{eqnarray*}
where $\tau^{n}_{i}$ are stopping times given by $\tau^{n}_{i}=\inf\{t>
t^{n}_i\dvtx Y^{n}_{t}> \bar Y^{n-1}_{t}\}\wedge t^{n}_{i+1}$, for
$i=0,\ldots, 2^{n}-1$.
By construction $((\bar Y^n,\bar Z^n))$ fulfills both \eqref
{eqstep1nprocess} and \eqref{eqthmmexynmono}, and $((\bar
Y^{n},\bar Z^{n}))\subset\mathcal{A}(\xi,g)$ with Lemma~\ref{lemadmpasting}.

\textit{Step} 3. \textit{Bound on $\int Z^n \,dW$.}
We now take the sequence $((Y^{n}, Z^{n}))$ fulfilling~\eqref
{eqstep1nprocess} and \eqref{eqthmmexynmono} and provide an
inequality which will enable us to use compactness arguments for $
( Z^n  )$ later in the proof.
More precisely, we argue that, for all $n\in{\mathbb N}$, it holds that
%
%e4.3 #&#
%
\begin{equation}
\label{eqthmmexistence1} \biggl\vert\int_{0}^{t}Z_{s}^n
\,dW_{s}\biggr\vert \leq B_t^n:=\bigl\vert
Y_t^1\bigr\vert +E \bigl[ \xi ^-| \mathcal{F}_t
\bigr]+E \bigl[ \xi^- \bigr]+\bigl\vert Y_0^1
\bigr\vert+A_t^n
\end{equation}
for all $t\in [ 0,T  ]$, where $ A^n_t$ is the positive
increasing process defined in Lemma~\ref{remdecompY}.
Moreover, it holds that
\[
E\bigl[A^n_T\bigr]\le Y^1_{0}-E[
\xi].
\]
Indeed, by the same arguments as in the proof of Lemma~\ref
{propproperties1}(2), recall $Y_0^n\leq Y_{0}^1$,
it follows that
%
%e4.4 #&#
%
\begin{equation}
\label{eqinequalityblabla01} \int_0^t
Z^n_s \,dW_s\ge-E \bigl[\xi^-|
\mathcal{F}_t \bigr]-Y_{0}^1.
\end{equation}
On the other hand, from $Y_{t}^n\leq Y_{t}^1$ and $-Y_{0}^n\leq E[\xi
^-]$, recall Lemma~\ref{remdecompY}, it follows that
%
%e4.5 #&#
%
\begin{equation}
\label{eqinequalityblabla02} \int_{0}^{t}Z_s^n
\,dW_s\leq Y_t^1+A_t^n-Y_0^n
\leq Y_t^1+A_t^n+E\bigl[\xi^-
\bigr].
\end{equation}
Combining \eqref{eqinequalityblabla01} and \eqref
{eqinequalityblabla02} yields \eqref{eqthmmexistence1}.
The $L^1$ bound on $A^n$ follows from $Y_{0}^n-A_{T}^n +\int_0^T Z^n_s
\,dW_s= \xi$, $Y_{0}^1\geq Y_{0}^n$ and the supermartingal property of
$\int Z^n \,dW$.

Note that if $(B^{n,\ast}_T)$ in \eqref{eqthmmexistence1} were
bounded in $L^1$, then, by means of the Burkholder--Davis--Gundy
inequality, $(Z^n)$ would be a bounded sequence in $\mathcal{L}^1$, and
we could apply \cite{delbaen03}, Theorem A, to find a sequence in the
asymptotic convex hull of $(Z^n)$ converging in $\mathcal{L}^1$ and
$P\otimes dt$-almost surely along some localizing sequence of stopping
times to some limit $Z \in\mathcal{L}^1$.
Here, even if $(A_T^{n,\ast})=(A_T^n)$ is uniformly bounded, this is,
however, not necessarily the case for $Y^{1,\ast}_T$ and $(E[\xi^-|
\mathcal{F}_{\cdot}])^*_{T}$, and this is the reason why we introduce
the following localization.

\textit{Step} 4. \textit{First localization.}
Due to our Brownian setting and since $\xi^-\in L^1$, we know that the
martingale $E[\xi^-| \mathcal{F}_{\cdot}]$, has a continuous
version; see \cite{Revuz1999}, Theorem~V.3.5.
Moreover, $Y^1$ is a c\`adl\`ag supermartingale and thus we may take
the localizing sequence
%
%e4.6 #&#
%
\begin{equation}
\label{eqsigmaseq} \sigma_{k}=\inf\bigl\{ t> 0 \dvtx\bigl \vert
Y_{t}^1\bigr\vert+ E \bigl[\xi^-| \mathcal {F}_t
\bigr]> k \bigr\} \wedge T,\qquad k \in\N,
\end{equation}
which is independent of $n \in\N$.
For a fixed $k \in\N$, inequality \eqref{eqthmmexistence1} yields
%
%e4.7 #&#
%
\begin{equation}
\label{eqproofexistenceboundZ} \biggl(\int Z^n1_{[0,\sigma_{k}]}\,dW
\biggr)^\ast_{T}\leq B^{k,n}\qquad\mbox{for all } n\in
\N,
\end{equation}
where $B^{k,n} = \vert Y^1_0\vert+E[\xi^-]+k+A^n_T$.
Due to $E[A^n_{T}]\le Y^1_{0}-E[\xi]$ we have
%
%e4.8 #&#
%
\begin{equation}
\label{eqBknL1bounded} \sup_{n\in{\mathbb N}} E \bigl[ B^{k,n}
\bigr]<\infty.
\end{equation}
Since $(B^{k,n})_{n\in\N}$ is a sequence of positive random variables,
we may apply \cite{DS94}, Lemma A1.1.
It provides a sequence $(\tilde B^{k,n})_{n\in\N}$ in the asymptotic
convex hull of $(B^{k,n})_{n\in\N}$, which converges almost surely to a
random variable $\tilde B^k\geq0$.
The $\tilde B^{k,n}$ inherit the integrability of the $B^{k,n}$, and
we can conclude with Fatou's lemma that
%
%e4.9 #&#
%
\begin{equation}
\label{eqtildeBintbar} E\bigl[\tilde B^k\bigr]< \infty.
\end{equation}
Let $\tilde Z^{k,n}$ be the convex combination of $(Z^n)$ corresponding
to $\tilde B^{k,n}$ so that
%
%e4.10 #&#
%
\begin{equation}
\label{eqthmmexistence3} \biggl(\int\tilde Z^{k,n}1_{[0,\sigma_{k}]}\,dW
\biggr)^\ast_{T}\leq \tilde B^{k,n} \qquad\mbox{for all
} n\in\N.
\end{equation}
%]

\textit{Step} 5. \textit{Second localization.}
The next two steps follow some known compactness arguments, which, in
the case of $\mathcal{L}^1$, can be found in \cite{delbaen03}.
For the sake of completeness we develop the argumentation.
Given an $m\in\N$, we start by taking a fast subsequence $(\tilde
B^{k,m,n})_{n\in\N}$ of $(\tilde{B}^{k,n})_{n\in\N}$ converging in
probability to $\tilde B^k$.
More precisely, we choose $(\tilde B^{k,m,n})_{n\in\N}$ such that
%
%e4.11 #&#
%
\begin{equation}
\label{eqfastseqm} P \bigl[\bigl\vert\tilde B^{k,m,n}-\tilde B^k
\bigr\vert\geq1 \bigr]\leq\frac
{2^{-n}}{m}.
\end{equation}
Consider now the stopping time $\tau^{k,m}$ given by
%
%e4.12 #&#
%
\begin{equation}
\label{eqtaum} \tau^{k,m} = \inf\biggl\{ t\geq0\dvtx \biggl(\int\tilde
Z^{k,m,n}1_{[0,\sigma _{k}]}\,dW \biggr)^*_{t}\geq m, \mbox{ for
some }n\in\mathbb{N} \biggr\} \wedge T,\hspace*{-35pt}
\end{equation}
where the sequence $(\tilde Z^{k,m,n}1_{[0,\sigma_{k}]})_{n\in\N}$ is
the subsequence of $(\tilde Z^{k,n}1_{[0,\sigma_{k}]})_{n\in\N}$
corresponding to the fast subsequence $(\tilde B^{k,m,n})_{n\in\N}$.
The definition of $\tau^{k,m}$ as well as the Burkholder--Davis--Gundy
inequality imply that the sequence of processes $(\tilde
Z^{k,m,n}1_{[0,\sigma_{k}]}1_{[0,\tau^{k,m}]})_{n\in\N}$ is bounded in
$\mathcal{L}^2$.
The Alaoglu--Bourbaki theorem and the Eberlein--\v{S}mulian theorem in
the Banach space $\mathcal{L}^2$\break imply the existence of $\hat
Z^{k,m}\in\mathcal{L}^2$, such that, up to a subsequence,\break $(\tilde
Z^{k,m,n}1_{[0,\sigma_{k}]}1_{[0,\tau^{k,m}]})_{n\in\N}$ converges
weakly to $\hat Z^{k,m}$.
As a consequence of the Hahn--Banach theorem, there exists a sequence
in the asymptotic convex hull of $(\tilde Z^{k,m,n}1_{[0,\sigma
_{k}]}1_{[0,\tau^{k,m}]})_{n\in\N}$, again denoted with $(\tilde
Z^{k,m,n}1_{[0,\sigma_{k}]}1_{[0,\tau^{k,m}]})_{n\in\N}$, which
converges in $\mathcal{L}^2$ to $\hat Z^{k,m}$.
By taking another subsequence we also have the $P\otimes dt$-almost
sure convergence.

\textit{Step} 6. \textit{$(\tau^{k,m})_{m\in\N}$ is a localizing
sequence of stopping times.}
We estimate as follows:
\begin{eqnarray*}
P \bigl[\tau^{k,m}=T \bigr]&=& P \biggl[ \biggl(\int\tilde
Z^{k,m,n}1_{[0,\sigma_{k}]}\,dW \biggr)^\ast_{T}<m,\mbox{ for all } n \in \N \biggr]
\\
&\geq&1-P \bigl[\tilde B^{k,m,n}\geq m, \mbox{ for some }n\in\N \bigr]
\\
&\geq&1- P \bigl[\bigl\{ \bigl\vert\tilde B^{k,m,n}-\tilde B^k
\bigr\vert\geq1, \mbox{ for some }n\in\N \bigr\} \cup\bigl\{ \tilde B^k >
m-1 \bigr\} \bigr]
\\
&\geq&1-\sum_{n}P \bigl[\bigl\llvert \tilde
B^{k,m,n}-\tilde B^k\bigr\rrvert \geq 1 \bigr]-P \bigl[
\tilde B^k > m-1 \bigr]
\\
&\geq&1-\frac1m-\frac{E [\tilde B^k +1 ]}{m} \mathop{\longrightarrow }_{m\rightarrow\infty} 1,
\end{eqnarray*}
where we used \eqref{eqthmmexistence3} in the second line and \eqref
{eqfastseqm}, the Markov inequality and the fact that $E[\tilde
B^k]<\infty$ in the last one.

\textit{Step} 7. \textit{Construction of the candidate $\hat Z$.}
For given $k,m>0$, we constructed in step 5 the process $\hat
Z^{k,m}$ as the $\mathcal{L}^2$ and $P\otimes dt$-almost sure limit of
a sequence in the asymptotic convex hull of $(\tilde
Z^{k,m,n}1_{[0,\sigma_{k}]}1_{[0,\tau^{k,m}]})_{n\in\N}$.
With $(\tilde B^{k,m,n})_{n\in\N}$ we denote the corresponding
subsequence of convex combinations of $(\tilde B^{k,m,n})_{n\in\N}$ and
note that $(\int\tilde Z^{k,m,n}1_{[0,\sigma_{k}]}\,dW)^*_{T}\leq
\tilde
B^{k,m,n}$, for all $n\in\N$, as in \eqref{eqthmmexistence3}.
Hence, by the same procedure as in step 5, we can find, for
$m'>m$, a fast subsequence $(\tilde Z^{k,m',n}1_{[0,\sigma
_{k}]})_{n\in
\N}$ in the asymptotic convex hull of $(\tilde Z^{k,m,n}1_{[0,\sigma
_{k}]})_{n\in\N}$ and a $\hat Z^{k,m'}\in\mathcal{L}^2$ such that
$(\tilde Z^{k,m',n}1_{[0,\sigma_{k}]}1_{[0,\tau^{k,m'}]})_{n\in\N}$
converges in $\mathcal{L}^2$ and $P\otimes dt$-almost surely to $\hat
Z^{k,m'}$.
We iterate this procedure and define $(\tilde Z^{k,n})_{n\in\N}$ as
the diagonal sequence $\tilde Z^{k,n} = \tilde Z^{k,n,n}$ and $\hat
Z^k$ as
%
%e4.13 #&#
%
\begin{equation}
\label{eqdeftildeZk} Z^k_{0}=0, \qquad \hat Z^k =
\sum_{m=1}^{\infty}\hat Z^{k,m}1_{]\tau
^{k,m-1},\tau^{k,m}]}.
\end{equation}
From $\hat Z^{k,m'}1_{[0,\tau^{k,m}]}=\hat Z^{k,m}$, for $m'>m$,
follows that $(\tilde Z^{k,n}1_{[0,\sigma_{k}]}1_{[0,\tau
^{k,n}]})_{n\in
\N}$ converges in $\mathcal{L}^2$ and $P\otimes dt$-almost surely to
$\hat Z^{k}$.
With the sequence $(\tilde Z^{k,n})_{n\in\N}$ and the process $\hat
Z^k$ at hand, we now diagonalize our program above with respect to $k$
and $n$.
As before, we get a diagonal sequence $\tilde Z^n=\tilde Z^{n,n}$, and
a process $\hat Z$ given by
%
%e4.14 #&#
%
\begin{equation}
\label{eqdefZ} \hat Z_{0} =0,\qquad \hat Z =\sum
_{k=1}^{\infty}1_{]\sigma_{k-1},\sigma
_{k}]}\hat Z^{k},
\end{equation}
such that
%
%e4.15 #&#
%
\begin{equation}
\label{eqconvdiagseqfinal} \tilde Z^{n}1_{[0,\tau_{n}]}
\mathop{\onrarrow{P\otimes dt\mathrm{\mbox{-}almost\ surely}}}_
{n\rightarrow\infty}\hat Z
\end{equation}
for $\tau_{n}=\sigma_{n}\wedge\tau^{n,n}$, where $\sigma_{n}$ and
$\tau
^{n,n}$ are as in \eqref{eqsigmaseq} and \eqref{eqtaum}, respectively.
For later reference, note that by construction it holds that $\hat
Z^{k',m}1_{[0,\sigma_{k}]}1_{[0,\tau^{k,m}]}=\hat Z^{k,m}$, as soon as
$ k'\geq k$ and also $\hat Z1_{[0,\sigma_{k}]}1_{[0,\tau^{k,m}]}=\hat
Z^{k,m}$.
Likewise $(\tilde Z^{n}1_{[0,\tau_{l}]})_{n\in\N}$ converges in
$\mathcal{L}^2$ and $P\otimes dt$-almost surely to $\hat Z^{l,l}$.
This yields, via the Burkholder--Davis--Gundy inequality, up to a subsequence,
%
%e4.16 #&#
%
\begin{eqnarray}
\label{eqconvstochint} \int_0^{t \wedge\tau_l}
\tilde{Z}^n_s \,dW_s \mathop{\longrightarrow}
_{n\rightarrow
+\infty}\int_0^{t \wedge\tau_l} \hat{Z}_s
\,dW_s
\nonumber
\\[-8pt]
\\[-8pt]
 \eqntext{\mbox{for all } t \in[0,T], P\mbox{-almost surely.}}
\end{eqnarray}
Hence, diagonalizing yields
%
%e4.17 #&#
%
\begin{equation}
\label{eqconvstochint2}\qquad \int_0^{t}
\tilde{Z}^n_s \,dW_s \mathop{\longrightarrow}_
{n\rightarrow+\infty}
\int_0^{t } \hat{Z}_s
\,dW_s\qquad \mbox{for all } t \in[0,T], P\mbox {-almost surely.}
\end{equation}

\textit{Step} 8. \textit{Monotone convergence to $\mathcal
{E}^{g}(\xi)$.}
Let $\tilde Y_{t}=\lim_{n}Y^n_t$, for $t\in[0,T]$, be the pointwise
monotone limit of the sequence $(Y^n)$.
By monotone convergence $\tilde Y$ is a supermartingale and, since our
filtration is right-continuous, by standard arguments we may define the
c\`adl\`ag supermartingale $\hat Y$ by setting $\hat{Y}_t=\lim_{s\downarrow t,s\in\mathbb{Q}}\tilde{Y}_s$, for all $t\in[0,T)$, and
$\hat Y_{T}=\xi$.
By construction, $\tilde{Y}_{t^i_n}=\hat{\mathcal{E}}^g_{t^i_n} (
\xi )$.
Hence, $\hat Y_{t}=\mathcal{E}^{g}_{t}(\xi)$, for all $t\in[0,T]$, and
%
%e4.18 #&#
%
\begin{equation}
\label{eqthmmexthatY} Y^n_t\geq\tilde{Y}_t
\geq\hat{\mathcal{E}}^g_{t} ( \xi )\geq
\mathcal{E}^{g}_{t}(\xi)\geq E [ \xi| \mathcal{F}_t
],
\end{equation}
where the third inequality follows from Proposition~\ref
{propsinnlosphisubmartingale}.
Now, the process $\mathcal{E}^{g}(\xi)$ is the natural candidate for
the value process of the minimal supersolution for two reasons.
It is c\`adl\`ag and it is dominated by $\hat{\mathcal{E}}^{g}(\xi)$
as \eqref{eqthmmexthatY} shows.
However, it is not clear a priori that the sequence $(Y^n)$ converges
to $\mathcal{E}^{g}(\xi)$ in some suitable sense.
Taking into account the additional structure provided by the
supermartingale property of the $Y^n$ we can prove nonetheless
%
%e4.19 #&#
%
\begin{equation}
\label{eqhatYmonlimit2} \mathcal{E}^{g}(\xi)=\hat Y=\lim
_{n\rightarrow\infty} Y^n,\qquad P\otimes dt\mbox{-almost surely.}
\end{equation}
To see this note first that by right continuity the limit $\tilde
Y_t=\lim_{n}Y^n_t$ is defined, for all $t\in[0,T]$, $P$-almost surely.
We now consider the sequence $((\tilde Y^n,\tilde Z^n))$ in the
asymptotic convex hull of $(Y^n, Z^n)$, which corresponds to the
sequence $(\tilde Z^n)$ constructed in step 7.
From the decomposition of the $Y^n$ (see Lemma~\ref{remdecompY}), we
obtain that $\tilde Y_{t}^n= \tilde Y_{0}^n- \tilde A_{t}^n+ \tilde
M^n_{t}$, for all $t\in[0,T]$.
Since $(\tilde Y_{t}^n)$ and $(\tilde M^n_{t})$ converge for all $t\in
[0,T]$, $P$-almost surely, the sequence $(\tilde A_{t}^n)$ also
converges, that is, there exists an increasing positive integrable
process $\tilde A$, such that $\lim_{n\to\infty} \tilde A^{{n}}_t =
\tilde A_t$, for all $t\in[0,T]$, $P$-almost surely.
Thus $\tilde Y_{t}=\tilde Y_{0}-\tilde A_{t}+\tilde M_{t}$, for all
$t\in[0,T]$.
Consequently, the jumps of $\tilde Y$ are given by the countably many
jumps of the increasing process $\tilde A$, which implies
\[
\hat Y_{t}=\tilde Y_{0}-\lim_{s\downarrow t,s\in\mathbb{Q}}
\tilde A_{s}+\tilde M_{t}\qquad \mbox{for all } t\in[0,T),\qquad \hat
Y_{T}=\xi.
\]
Moreover, the jump times of the c\`adl\`ag process $\hat Y$ are
exhausted by a sequence of stopping times $(\rho^j)\subset\mathcal{T}$,
which coincide with the jump times of~$\tilde A$.
Therefore, $\hat Y=\tilde Y$, $P\otimes dt$-almost surely, which
implies \eqref{eqhatYmonlimit2}.

\textit{Step} 9. \textit{Verification.}
Let us now show that $(\mathcal{E}^{g}(\xi),\hat{Z})\in\mathcal
{A} ( \xi,g  )$, which, by means of \eqref{eqthmmexthatY},
would end the proof.
We start with the verification of \eqref{eqcentralineq} under the
assumption~\textsc{(Inc)}.
Due to \eqref{eqhatYmonlimit2} there exists a set $B\subset\Omega
\times[0,T]$ with $P\otimes dt(B^c)=0$, such that $\mathcal
{E}^{g}_t(\xi)(\omega)=\lim_{n\rightarrow\infty} Y^n_t(\omega)$, for
all $(\omega,t)\in B$.
Hence, there exists a set $A\subset\{\omega\dvtx  (\omega,t)\in B, \mbox
{for some }t\}$, with $P(A)=1$, such that, for all $\omega\in A$, the
set $I(\omega)=\{t\in[0,T]\dvtx (\omega,t)\in B\}$ is a Lebesque set of
measure $T$ and $\mathcal{E}^{g}_t(\xi)(\omega)=\lim_{n\rightarrow
\infty} Y^n_t(\omega)$, for all $t\in I(\omega)$.
In the following we suppress the dependence of $I$ on $\omega$ and
just keep in mind that $s$ and $t$ may depend on $\omega$.
Let $s,t\in I$ with $s\leq t$.
% Due to \eqref{eqhatYmonlimit2} there exists a set $\mathcal{A}
% In the following we suppress the dependence of $\mathcal{I}$ on $
% Let $s,t\in\mathcal{I}^c$ with $s\leq t$.
By using \eqref{eqconvstochint2}, the $P\otimes dt$-almost sure
convergence of $\tilde Z^n1_{[0,\tau^n]}$ to $\hat Z$, and Fatou's
lemma we obtain
%
%e4.20 #&#
%e4.21 #&#
%
\begin{eqnarray}
\label{thmmexistencestep52} &&\mathcal{E}^{g}_{s}-\int
_{s}^{t}g_{u}\bigl(
\mathcal{E}^{g}_{u},\hat Z_u\bigr)\,du+\int
_{s}^{t}\hat Z_u \,dW_u
\nonumber
\\[-8pt]
\\[-8pt]
\nonumber
&&\qquad\geq\limsup_{n} \biggl( \tilde Y_{s}^n-
\int_{s}^{t}g_{u}\bigl(\mathcal
{E}^{g}_{u},\tilde Z_u^n1_{[0,\tau^n]}(u)
\bigr)\,du+\int_{s}^{t}\tilde Z_u^n
\,dW_u \biggr),
\end{eqnarray}
where $\tilde Y^{n}$ denotes the convex combination of $(Y^n)$
corresponding to $\tilde Z^{n}$.
We denote by $\lambda^{(n)}_{i}$, $n\leq i\leq M^{(n)}$, $\lambda
^{(n)}_{i}\geq0, \sum_{i}\lambda^{(n)}_{i}=1$ the convex weights of
$\tilde Z^n$.
Since our generator fullills \textsc{(Con)}, and since, for $n$ large
enough, we have $\tilde Z_u^n1_{[0,\tau^n]}(u)=\tilde Z_u^n$, for all
$s\leq u\leq t$, we may further estimate the above by
\begin{eqnarray*}
&&\mathcal{E}^{g}_{s}-\int_{s}^{t}g_{u}
\bigl(\mathcal{E}^{g}_{u},\hat Z_u\bigr)\,du+\int
_{s}^{t}\hat Z_u \,dW_u
\\
&&\qquad\geq\limsup_{n}\sum_{i=n}^{M^{(n)}}
\lambda^{(n)}_{i} \biggl( Y_{s}^i-
\int_{s}^{t}g_{u}\bigl(
\mathcal{E}^{g}_{u},Z_u^i\bigr)\,du+
\int_{s}^{t}Z_u^i
\,dW_u \biggr).
\end{eqnarray*}
Since $Y^{i}_{t}\geq\hat{\mathcal{E}}^g_{t}(\xi)\geq\mathcal
{E}^{g}_{t}(\xi)$, for all $t\in[0,T]$, and $i\in\N$, we use \textsc{(Inc)} and the fact that the $(Y^n,Z^n)$ are supersolutions to conclude
%
%e4.22 #&#
%e4.23 #&#
%e4.24 #&#
%
\begin{eqnarray}
\label{thmmexistencestep54}&& \mathcal{E}^{g}_{s}-\int
_{s}^{t}g_{u}\bigl(
\mathcal{E}^{g}_{u},\hat Z_u\bigr)\,du+\int
_{s}^{t}\hat Z_u \,dW_u
\nonumber\\
&&\qquad\geq\limsup_{n}\sum_{i=n}^{M^{(n)}}
\lambda^{(n)}_{i} \biggl( Y_{s}^i-
\int_{s}^{t}g_{u}\bigl(Y_u^{i},Z_u^i
\bigr)\,du+\int_{s}^{t}Z_u^i
\,dW_u \biggr)
\\
&&\qquad\geq\limsup_{n}\sum_{i=n}^{M^{(n)}}
\lambda^{(n)}_{i} Y_{t}^i= \limsup
_{n}\tilde Y_{t}^n=\limsup
_{n} Y_{t}^n = \mathcal{E}^{g}_{t}.\nonumber
\end{eqnarray}
As for the case of $s,t \in I^c$, with $s \leq t$, we approximate them
both from the right with some sequences $(s^n)\subset I$ and
$(t^n)\subset I$, such that $s^n\downarrow s$, $t^n\downarrow t$,
$s^n\leq t^n$.
For each $s^n$ and $t^n$, \eqref{thmmexistencestep54} holds.
Passing to the limit\vspace*{1pt} by using the right-continuity of $\mathcal
{E}^{g}$ and the continuity of $-\int g(\mathcal{E}^{g},\hat Z)\,du+\int
\hat Z \,dW$, we deduce that \eqref{thmmexistencestep54} holds for all
$s,t \in [ 0,T  ]$ with $s\leq t$.

It remains to show admissibility of $\hat Z$.
By means of \eqref{thmmexistencestep54}, \eqref{eqthmmexthatY}
and positivity of $g$ it holds that
%
%e4.25 #&#
%
\begin{equation}
\label{eqexthmmadm} \int_0^{t}\hat
Z_s\,dW_s\geq E [\xi| \mathcal{F}_t ]-
\mathcal{E}_{0}.
\end{equation}
Being bounded from below by a martingale, the continuous local
martingale $\int_{}^{}\hat{Z}\,dW$ is by Fatou's lemma a supermartingale,
and thus $\hat{Z}$ is admissible.
Hence, the proof under assumptions \textsc{(Pos)}, \textsc{(Con)} and
\textsc{(Inc)} is complete.
%
% by using $\tilde Y_{t}\geq E[\xi\mid\mathcal{F}_{t}]$, recall that
%the $Y^n$ are supermartingales, instead of $\tilde Y_{t}\geq\hat
%Y_{t}$ in \eqref{thmmexistencestep54}.
% From positivity of $g$ and $\hat Y_{0}\leq\mathcal{E}^g_{0}(\xi)$
%follows $\int_0^{t}\hat Z_sdW_s\geq E\left[\xi\Mid\mathcal{F}_t\right]-
%admissible.
% Finally, $(\hat Y,\hat Z)\in\mathcal{A}(\xi,g)$ and in turn $\hat Y\ge

The proof under \textsc{(Dec)} replacing \textsc{(Inc)} only differs in
the verification of \eqref{eqcentralineq}.
Indeed, instead of only approximating $\hat Z$ in the Lebesgue
integral, we approximate $\mathcal{E}^{g}(\xi)$ $P\otimes dt$-almost
surely with the sequence $(Y^n)$ as well, that is, \eqref
{thmmexistencestep52} becomes, by means of \eqref
{eqhatYmonlimit2} and Fatou's lemma,
\begin{eqnarray*}
&&\mathcal{E}^{g}_{s}-\int_{s}^{t}g_{u}
\bigl(\mathcal{E}^{g}_{u},\hat Z_u\bigr)\,du+\int
_{s}^{t}\hat Z_u \,dW_u
\\
&&\qquad\geq\limsup_{n} \biggl( \tilde Y_{s}^n-
\int_{s}^{t}g_{u}\bigl(Y_u^n,
\tilde Z_u^n1_{[0,\tau^n]}(u)\bigr)\,du+\int
_{s}^{t}\tilde Z_u^n
\,dW_u \biggr).
\end{eqnarray*}
This entails, by monotonicity of the sequence $(Y^n)$ and the fact that
the convex combinations in $\tilde Z^n$ consist of elements of $(Z^i)$
with index greater or equal than~$n$, that we may write $-\int_{s}^{t}g_{u}(Y_u^n, Z_u^{i})\,du\geq-\int_{s}^{t}g_{u}( Y_u^{i},
Z_u^{i})\,du$ in \eqref{thmmexistencestep54}, and this ends the proof.
\end{pf*}

%re4.2 #&#
\begin{rem}\label{remarkvolatarpo}
Note that the existence theorem also holds if we additionally take
into account a volatility process in the stochastic integral.
More precisely, consider a progressively measurable process $\sigma
\dvtx \Omega\times [ 0,T  ]\to\mathbb{S}^{>0}_d$, where
$\mathbb
{S}^{>0}_d$ denotes the set of strictly positive definite $d\times d$
matrices and define $\mathcal{L}^{\sigma}$ as the set of progressively
measurable processes $Z \dvtx \Omega\times [ 0,T  ]\to\R
^{1\times d}$ such that $Z \sigma^{1/2} \in\mathcal{L}$.
Analogously to the previous setting, given a generator $g$ and a
terminal condition $\xi\in L^0$, we say that $(Y,Z)\in\mathcal
{S}\times\mathcal{L}^{\sigma}$ is a supersolution of the BSDE under
volatility~$\sigma$ if
%
%e4.26 #&#
%
\begin{equation}
Y_s-\int_{s}^{t}g_u (
Y_u,Z_u )\,du +\int_{s}^{t}Z_u
\sigma _u^{1/2} \,dW_u\geq Y_t\quad
\mbox{and} \quad Y_T\geq\xi \label{eqcentralineqh}
\end{equation}
for all $0\leq s \leq t\leq T$.
We say that the control process is admissible if\break $\int_{}^{}Z\sigma
^{1/2} \,dW$ is a supermartingale, and define
%
%e4.27 #&#
%
\begin{equation}\quad
\mathcal{A} ( \xi,g,\sigma )=\bigl\{ (Y,Z)\in\mathcal {S}\times
\mathcal{L}^\sigma\dvtx  Z\mbox{ is admissible and }\eqref {eqcentralineqh}
\mbox{ holds} \bigr\} \label{eqcentralAh}
\end{equation}
as well as
%
%e4.28 #&#
%
\begin{equation}
\hat{\mathcal{E}}^{g,\sigma}_{t} ( \xi )=\essinf\bigl\{
Y_t\dvtx (Y,Z)\in\mathcal{A} ( \xi,g,\sigma ) \bigr\},\qquad t \in [ 0,T ].
\label{}
\end{equation}
We can formulate the following existence theorem.
%
%th4.3 #&#
\begin{thmm}
Let $g$ be a generator fulfilling \textsc{(Pos)}, \textsc{(Lsc)}, \textsc{(Con)} and either \textsc{(Inc)} or \textsc{(Dec)} and $\xi\in L^0$ be a
terminal condition, such that \mbox{$\xi^-\in L^1$}.

If $\mathcal{A}(\xi,g,\sigma) \neq\varnothing$, then there exists a
unique minimal supersolution $(\hat Y,\hat Z)\in\mathcal{A}(\xi
,g,\sigma)$.
Moreover, $\mathcal{E}^{g,\sigma} ( \xi )$ is the value
process of the minimal supersolution, that is, $(\mathcal{E}^{g,\sigma
} ( \xi ),\hat{Z})\in\mathcal{A}(\xi,g,\sigma)$.
\end{thmm}
The proof follows exactly the same scheme as the proof of Theorem~\ref
{thmmexistence} with a compactness argument in the Hilbert space
$\mathcal{L}^{2,\sigma}$, the set of processes in $\mathcal
{L}^\sigma$
such that $E[\int_{0}^{T}(Z_u \sigma_u^{1/2})(Z_u \sigma
_u^{1/2})^\top
du]< +\infty$, instead of $\mathcal{L}^2$.
\end{rem}
Theorem~\ref{thmmexistence} ensures the existence and uniqueness of
the minimal supersolution which is c\`adl\`ag.
The following proposition provides a condition under which $\mathcal
{E}^g(\xi)$ is in fact continuous.
%
%pr4.4 #&#
\begin{prop}\label{propbsde02}
Let $g$ be a generator fulfilling \textsc{(Pos)}, \textsc{(Lsc)}, \textsc{(Con)} and either \textsc{(Inc)} or \textsc{(Dec)} and $\xi\in L^0$ be a
terminal condition, such that $\xi^-\in L^1$.
Suppose that $\mathcal{A}(\xi,g) \neq\varnothing$.
Assume that for any $\zeta\in L^\infty ( \mathcal{F}_{\tau}
)$, $\tau\in\mathcal{T}$, there exist $Y \in\mathcal{S}$ and an
admissible $Z\in\mathcal{L}$, which solve the backward stochastic
differential equation
\[
Y_t-\int_{t}^{\tau}g_{s} (
Y_s,Z_s )\,ds+\int_{t}^{\tau}Z_s
\,dW_s=\zeta\qquad \mbox{for all } t \in[0,\tau].
\]
Then $\mathcal{E}^g ( \xi )$ is continuous.
\end{prop}
\begin{pf}
In view of Theorem~\ref{thmmexistence}, there exists $\hat{Z}\in
\mathcal{L}$ such that $(\mathcal{E}^g, \hat{Z})\in\mathcal{A}(\xi,g)$.
Hence, $\mathcal{E}^g$ can only have downward jumps.
Assume that $\mathcal{E}^g$ has a negative jump, that is, $P[\tau\leq
T]>0$, for the stopping time $\tau=\inf\{t>0 \dvtx \Delta\mathcal
{E}^g_t<0\}$.
We then fix $m$ big enough such that the stopping time $\tau^m=\inf\{
t>0\dvtx |\mathcal{E}^g_t|>m\}\wedge\tau$ satisfies $P[ \{-m<\Delta
\mathcal{E}^g_{\tau^m}<0\}\cap\{\tau^m=\tau\}]>0$.
Since $\mathcal{E}^g$ is continuous on $[ 0,\tau[$, and $\mathcal
{E}^g$ has only negative jumps, $\mathcal{E}^g_{\tau^m}\vee-m \in
L^\infty ( \mathcal{F}_{\tau^m}  )$.
By assumption there exist $\bar Y \in\mathcal{S}$ and an admissible
$\bar Z\in\mathcal{L}$ such that
\[
\bar{Y}_s+\int_{s}^{\tau^m}g_{u}
( \bar{Y}_{u},\bar{Z}_u )-\int_{s}^{\tau^m}
\bar{Z}_u \,dW_u=\mathcal{E}^g_{\tau^m}
\vee-m\qquad \mbox{for all } s\in\bigl[0,\tau^m\bigr].
\]
Similarly to Lemma~\ref{lemadmpasting}, we derive $(\bar{Y}1_{[0,\tau
^m[}+\mathcal{E}^g1_{ [ \tau^m,T  ]},\bar{Z}1_{ [
0,\tau^m
]}+\hat Z1_{ ]\tau^m,T  ]})\in\mathcal{A} (
\xi,g
)$.
Hence, by optimality of $\mathcal{E}^g$ in $\mathcal{A}(\xi,g)$, it
holds that $\mathcal{E}^g\leq\bar{Y}1_{[0,\tau^m[}+\mathcal
{E}^g1_{ [ \tau^m,T  ]}$.
Moreover, we have
\[
\mathcal{E}^g_{\tau^m-}>\mathcal{E}^g_{\tau^m}
\vee-m=\bar {Y}_{\tau
^m}=\bar{Y}_{\tau^m-} \qquad\!\!\mbox{on the set }\bigl\{
-m<\Delta\mathcal {E}^g_{\tau^m}<0 \bigr\} \cap\bigl\{
\tau^m=\tau \bigr\}.
\]
Hence, for the stopping time $\hat\tau=\inf\{t>0 \dvtx \mathcal
{E}^g_t>\bar
{Y}_t\}\wedge\tau^m$ we deduce $P[\hat\tau<\tau^m]>0$, since the
processes $\mathcal{E}^g$ and $\bar Y$ are continuous on $[0,\tau^m[$.
But then $\mathcal{E}^g\nleq\bar Y$ on $[0,\tau^m[$, which is a
contradiction.
\end{pf}

Under the assumptions of Theorem~\ref{thmmexistence}, $\mathcal{E}^g$
is the value process of the minimal supersolution with a control
process $\hat{Z}$ in $\mathcal{L}$ which defines a supermartingale.
Next we address the following question: under which conditions does the
control process have enough integrability in order to define a true
martingale? That is, when does $\hat{Z}$ belong to some $\mathcal
{L}^p$, for $p\geq1$?
Defining
%
%e4.29 #&#
%
\begin{equation}
\mathcal{A}^p ( \xi,g ):=\bigl\{ ( Y,Z ) \in \mathcal {A} ( \xi,g )\dvtx Z \in\mathcal{L}^p \bigr\}, \label{setAp}
\end{equation}
this means that $(\mathcal{E}^g ( \xi ),\hat{Z}) \in
\mathcal
{A}^p ( \xi,g  )$.
Peng \cite{peng99} provides a positive answer to this question in the case
where $p=2$, the terminal condition $\xi\in L^2$ and the generator is
not necessarily positive but Lipschitz.
Compare this also with Cheridito and Stadje \cite{cs11} for
supersolutions of BSDEs where
the control process is in BMO, if the terminal condition is a bounded
lower semicontinuous function of the Brownian motion and the generator
is convex in $z$ and Lipschitz and increasing in $y$.
Here, we provide an answer to the case where $p=1$ in the context of
Section~\ref{sectarpointro}.
Given a terminal condition $\xi$, obtaining $\mathcal{E}^g ( \xi
)$ as a minimal solution with a control process within $\mathcal
{L}^1$ comes at two costs.
Indeed, a stronger integrability condition on the terminal value is
required; that is, we impose that $ (E[\xi^-| \mathcal{F}_{\cdot
}] )^\ast_{T}\in L^1$.
As for the second cost, $\mathcal{A}^1 ( \xi,g  )\neq
\varnothing
$ is also required, which, in view of $\mathcal{A}^1 ( \xi,g
)\subset\mathcal{A} ( \xi,g  )$, is also a stronger assumption.

%th4.5 #&#
\begin{thmm}\label{thmmexistenceH1}
Suppose that the generator $g$ fulfills \textsc{(Pos)}, \textsc{(Lsc)},
\textsc{(Con)} and either \textsc{(Inc)} or \textsc{(Dec)}.
Let $\xi\in L^0$ be a terminal condition, such that $ (E[\xi
^-| \mathcal{F}_{\cdot}] )^\ast_{T} \in L^1$.
If $\mathcal{A}^1 (\xi,g ) \neq\varnothing$, then there exists
a unique minimal supersolution $(\hat Y,\hat Z)\in\mathcal{A}^1(\xi,g)$.
Moreover, $\mathcal{E}^g ( \xi )$ is the value process of
the minimal supersolution, that is, $(\mathcal{E}^g ( \xi
),\hat{Z})\in\mathcal{A}^1(\xi,g)$.
\end{thmm}
%
%re4.6 #&#
\begin{rem}\label{remsupYintbar}
As in Section~\ref{sectarpointro}, note that for $(Y,Z)\in\mathcal
{A}^1 (\xi,g )$, the value process $Y$ is a supermartingale
with terminal value greater or equal than $\xi$.
Moreover, we have $Y^\ast_T\in L^1$.
Indeed, by using the decomposition \eqref{eqdecompY}, we derive
$Y^\ast_{t}\leq\vert Y_{0}\vert+A_{T}+ (\int Z \,dW )^\ast_T$.
We further have $A_{T}\leq Y_{0}+\int_{0}^TZ_{s}\,dW_{s}-\xi$ and thus
$E [\vert A_{T}\vert ]\leq Y_{0}+E [\xi^- ]$.
Consequently,
\[
E \bigl[Y^\ast_T \bigr]\leq\vert Y_{0}\vert+E
\bigl[\xi^- \bigr]+Y_{0}+E \biggl[ \biggl( \int_{}^{}Z
\,dW \biggr)^\ast_T \biggr]<\infty.
\]
\end{rem}
\begin{pf*}{Proof of Theorem~\ref{thmmexistenceH1}}
Since $\mathcal{A}^1 (\xi,g  )\subset\mathcal{A}(\xi,g)$,
the assumption\break $\mathcal{A}^1 ( \xi,g  )\neq\varnothing$
implies the existence of $\hat{Z}\in\mathcal{L}$ such that
$(\mathcal
{E}^g(\xi),\hat{Z})\in\mathcal{A} ( \xi,g  )$.
We are left to show that $\hat{Z}\in\mathcal{L}^1$.
Since $\mathcal{A}^1 ( \xi,g  )\neq\varnothing$, we can
suppose in the proof of Theorem~\ref{thmmexistence} that $(Y^1,Z^1)\in
\mathcal{A}^1 ( \xi,g  )$.
Since \eqref{eqthmmexistence1} holds for $(\mathcal{E}^g ( \xi
),\hat{Z})$, instead of $(Y^n,Z^n)$, we have
%
%e4.30 #&#
%
\begin{equation}
\label{eqcentralinequH1} \biggl(\int\hat{Z}\,dW \biggr)^\ast_{T}
\leq\bigl\vert Y^1_0\bigr\vert+E \bigl[\xi ^- \bigr]+\hat
A_{T}+\bigl(Y^1\bigr)^\ast_{T}+
\bigl(E \bigl[\xi^-| \mathcal {F}_{\cdot
} \bigr] \bigr)^\ast_{T},
\end{equation}
where $0\leq E[ \hat A_{T} ] \leq E[ \xi]-Y_0^1$.
Since $(E[ \xi^-| \mathcal{F}_{\cdot} ])^\ast_{T}\in L^1$, by means
of Remark~\ref{remsupYintbar}, the right-hand side of \eqref
{eqcentralinequH1}, is in $L^1$.
Thus, by means of the Burkholder--Davis--Gundy inequality, $\hat{Z}$
belongs to $\mathcal{L}^1$.
\end{pf*}

%
%s4.2 #&#
\subsection{Stability results}\label{secstab}\label{subsecstabterminal}

In this section we address the stability of $\hat{\mathcal
{E}}^g(\cdot
)$ with respect to perturbations of the terminal condition or the generator.
First we show that the functional $\hat{\mathcal{E}}^g_{0}$ is not only
defined on the same domain as the usual expectation, but also shares
some of its main properties, such as Fatou's lemma as well as a
monotone convergence theorem.
%
%th4.7 #&#
\begin{thmm}\label{thmmmonconvfatou}
Suppose that the generator $g$ fulfills \textsc{(Pos)}, \textsc{(Lsc)},
\textsc{(Con)} and either \textsc{(Inc)} or \textsc{(Dec)}.
Let $(\xi^n)$ be a sequence in $L^0$, such that $\xi^n\geq\eta$, for
all $n\in\N$, where $\eta\in L^1$.
\begin{itemize}%[ font= \normalfont\itshape]
\item Monotone convergence: if $(\xi^n)$ is
increasing $P$-almost surely to $\xi\in L^0$, then $\hat{\mathcal
{E}}^g_{0}(\xi)=\lim_{n}\hat{\mathcal{E}}^g_{0}(\xi^n)$.
\item Fatou's lemma: $\hat{\mathcal{E}}^g_{0}(\liminf_{n} \xi
^n)\leq\liminf_{n}\hat{\mathcal{E}}^g_{0}(\xi^n)$.
\end{itemize}
\end{thmm}
\begin{pf}
\textit{Monotone convergence}:
From Proposition~\ref{propproperties1} and by monotonicity, it
follows that $\hat{\mathcal{E}}^g(\xi^n)\leq\hat{\mathcal
{E}}^g(\xi
^{n+1})\leq\cdots\leq\hat{\mathcal{E}}^g(\xi)$.
Hence, we may define $\hat Y_0=\lim_{n}\hat{\mathcal{E}}^g_0(\xi^n)$.
Note that $\hat Y_0\leq\hat{\mathcal{E}}^g_0(\xi)$.
If $\hat Y_{0}=+\infty$, then also $\hat{\mathcal{E}}^g_{0}(\xi
)=+\infty$, and there is nothing to prove.
Suppose now that $\hat Y_{0}<\infty$.
This implies that $\mathcal{A}(\xi^n,g)\neq\varnothing$, for all
$n\in\N$.
Since $\xi^n\geq\eta$, Proposition~\ref
{propsinnlosphisubmartingale} yields $(\xi^n) \subset L^1$ and
$(\mathcal{E}^g(\xi^n))$ is a well-defined\vspace*{1pt} increasing sequence of c\`
adl\`ag supermartingales.
We define $Y_t=\lim_{n}\mathcal{E}^g_t(\xi^n)$, for all $t\in[0,T]$.
Note that $Y_0=\hat Y_0$.
We show that $Y$ is a c\`adl\`ag supermartingale.

To this end, note that the sequence $ (\mathcal{E}^g(\xi
^n)-\mathcal{E}^g(\xi^1) )$ is positive and increases to
$Y-\mathcal{E}^g(\xi^1)$. Therefore monotone convergence yields
\[
0\leq E\bigl[Y_t-\mathcal{E}^g_t\bigl(
\xi^1\bigr)\bigr]=\lim_n E\bigl[
\mathcal{E}^g_t\bigl(\xi ^n\bigr)-
\mathcal{E}^g_t\bigl(\xi^1\bigr)\bigr].
\]
The supermartingale property of $\mathcal{E}^g(\xi^n)$ implies that
$E[\mathcal{E}^g_t(\xi^n)]\leq\mathcal{E}^g_0(\xi^n)\leq Y_0$.
Furthermore, $E[\xi^1]\leq E[\mathcal{E}^g_t(\xi^1)]\leq Y_0$ and thus
\[
0\leq E\bigl[Y_t-\mathcal{E}^g_t\bigl(
\xi^1\bigr)\bigr]\leq-E\bigl[\xi^1\bigr]+Y_0<+
\infty.
\]
From $\mathcal{E}^g_t(\xi^1) \in L^1$, we deduce that $Y_t \in L^1$.
Since $\xi= Y_T$, this implies in particular that $\xi\in L^1$.
The supermartingale property follows by a similar argument.
% Indeed, again by monotone convergence we have
% \begin{eqnarray*}
% E\left[ \tilde Y_{t}\Mid\mathcal{F}_{s}\right]=E\left[\tilde Y_{t}-
% = \lim_{n}E\left[\mathcal{E}^g_{t}(\xi^n)- \mathcal{E}^g_{t}(\xi^1)
% \leq\lim_{n}\mathcal{E}^g_{s}(\xi^n)=\tilde Y_{s}.
% \end{eqnarray*}
Moreover, \cite{Dellacherie1982}, Theorem VI.18 implies that $Y$ is
indistinguishable from a c\`adl\`ag process.
Hence, $Y$ is a c\`adl\`ag supermartingale.

Theorem~\ref{thmmexistence} provides a sequence of optimal controls
$(Z^n)$ such that $(\mathcal{E}^g(\xi^n),Z^n) \in\mathcal{A}(\xi
^n,g)$, for all $n\in\N$.
Now we apply the procedure introduced in the proof of Theorem~\ref
{thmmexistence} and obtain a candidate control process $\hat Z$.
The only notable difference in the proof, except for the fact that $Y$
is already c\`adl\`ag, is that, here, the sequence $(\mathcal{E}^g(\xi
^n))$ is increasing instead of decreasing.
Thus, the c\`adl\`ag supermartingales $Y$ and $\mathcal{E}^g(\xi^1)$
serve as upper and lower bounds, respectively.
Consequently, we replace $Y^1$ by $Y$ and $E[\xi^-| \mathcal
{F}_{\cdot}]$ by $\mathcal{E}^g(\xi^1)$ in the key inequality~\eqref
{eqthmmexistence1}.
The verification follows exactly the same argumentation as in the
proof of Theorem~\ref{thmmexistence} for both monotonicity assumptions
\textsc{(Inc)} and \textsc{(Dec)}.
%The case \ref{cond03} is treated similarly, that is by jointly
%approximating $(\hat Y, \hat Z)$ with $((Y^{M^{(n)}},\tilde Z^n))$,
%where $M^{(n)}$ is the maximum of the indices of the convex weights in
%the representation of $\tilde Z^n$.
Finally, to get the admissibility of $\hat Z$, we denote with $(\tilde
\xi^n)$ the sequence of convex combinations of $(\xi^n)$ corresponding
to $(\tilde Z^n)$.
Monotonicity of the sequence $(\xi^n)$ implies $\xi^1\leq\tilde\xi
^n\leq\xi$, for all $n\in\N$.
We may and do switch to a subsequence such that $(\tilde\xi^n)$ is
increasing as well.
Now, fix an arbitrary $t\in[0,T]$.
Dominated convergence implies the $L^1$-convergence $\lim_{n}E[ \tilde
\xi^n| \mathcal{F}_{t}]=E[ \xi| \mathcal{F}_{t}]$.
Hence, we may select a subsequence such that we have $P$-almost sure
convergence.
Similarly to \eqref{eqexthmmadm} this implies
\[
Y_{0}-\int_{0}^{t}g_{u}(Y_u,
\hat Z_u)\,du+\int_{0}^{t}\hat
Z_u \,dW_u \geq \limsup_{n} E \bigl[
\tilde\xi^{n}| \mathcal{F}_{t} \bigr]= E [ \xi |
\mathcal{F}_{t} ].
\]
As before, this entails that $(Y, \hat Z)\in\mathcal{A}(\xi,g)$.
Hence, from $\mathcal{A}(\xi,g)\neq\varnothing$ and $\xi^-\in L^1$, we
derive by Theorem~\ref{thmmexistence} that there exists a control
process $Z$ such that $(\mathcal{E}^g(\xi),Z)\in\mathcal{A}(\xi,g)$.
In particular this yields $Y_{0}= \mathcal{E}^g_{0}(\xi)$, that is,
$\lim_{n}\mathcal{E}^g_{0}(\xi^n)= \mathcal{E}^g_{0}(\xi)$, since
otherwise $\mathcal{E}^g_{0}(\xi)$ were not optimal.

\textit{Fatou's lemma}:
The result follows by applying monotone convergence.\break
Indeed, denote by $\zeta^n$ the random variables $\zeta^n=\inf_{k\geq
n}\xi^k$.
Then from\break $\liminf_{n}\xi^n=\lim_{n}\zeta^n$, $\zeta^n\geq\eta$,
$\zeta^n\leq\xi^n$, for all $n\in\N$, and monotone convergence follows
\[
\hat{\mathcal{E}}^g_{0}\Bigl(\liminf_{n}
\xi^n\Bigr)=\hat{\mathcal {E}}^g_{0}\Bigl(\lim
_{n}\zeta^n\Bigr)=\lim_{n}
\hat{\mathcal{E}}^g_{0}\bigl(\zeta^n\bigr)\leq
\liminf_{n}\hat{\mathcal{E}}^g_{0}
\bigl(\xi^n\bigr).
\]
\upqed\end{pf}
%
%re4.8 #&#
\begin{rem}
An inspection of the proof of Theorem~\ref{thmmmonconvfatou} shows
that under the assumptions implying monotone convergence, if $\lim_{n}\hat{\mathcal{E}}^g_0 ( \xi^n  )< +\infty$, then
$\mathcal
{A} ( \xi,g  )\neq\varnothing$, and $\mathcal
{E}^g_{t} ( \xi
^n  )$ converges $P$-almost surely to $\mathcal{E}^g_t (
\xi
)$, for all $t \in [ 0,T  ]$.

Similarly, given a sequence $((Y^n,Z^n)) \subset\mathcal{A}(\xi,g)$
such that $(Y^n)$ is increasing and $\lim_{n}Y^n_{0}<\infty$, then
there exists a control process $Z\in\mathcal{L}$ such that $(Y,Z)\in
\mathcal{A}(\xi,g)$, where $Y_{t}$ is the $P$-almost sure limit of
$(Y^n_{t})$, for all $t\in[0,T]$.
\end{rem}
A consequence of the preceding theorem is the following result on
$L^1$-lower semicontinuity.
%
%th4.9 #&#
\begin{thmm}\label{thmmusc}
Let $g$ be a generator fulfilling \textsc{(Pos)}, \textsc{(Lsc)}, \textsc{(Con)} and either \textsc{(Inc)} or \textsc{(Dec)}.
Then $\hat{\mathcal{E}}^g_0$ is $L^1$-lower semicontinuous.
\end{thmm}
\begin{pf}
Let $(\xi^n)$ be a sequence of terminal conditions, which converges in
$L^1$ to a random variable $\xi$.
Suppose that there exists a subsequence $( \tilde{\xi}^n )\subset
(\xi
^n)$ such that $(\hat{\mathcal{E}}^g_{0}(\tilde{\xi}^n))$ converges to
some real $a< \hat{\mathcal{E}}^g_{0}(\xi)$.
We can assume, up to another fast subsequence, that $\|\tilde{\xi
}^n-\xi\|_{L^1}\leq2^{-n}$, for all $n\in\N$.
Consider now the sequence $(\zeta^n)$, with $\zeta^n$ given by
\[
\zeta^n=\xi-\sum_{k\geq n}\bigl(\tilde{
\xi}^k-\xi\bigr)^-.
\]
Clearly, $\zeta^n\in L^1$ and $\zeta^n\leq\zeta^{n+1}\leq\cdots
\leq
\xi$.
Moreover, $(\zeta^n)$ converges in $L^1$ to $\xi$, and, since it is
increasing, it converges also $P$-almost surely.
Thus, from Theorem~\ref{thmmmonconvfatou}, we get $\lim_{n}\hat
{\mathcal{E}}^g_{0}(\zeta^n)= \hat{\mathcal{E}}^g_{0}(\xi)$.
Now, $\zeta^n\leq\xi- (\tilde{\xi}^n-\xi)^-+(\tilde{\xi}^n-\xi
)^+\leq\tilde\xi^n$ and monotony of the functional $\hat{\mathcal
{E}}^g_{0}$ imply $a=\lim_{n}\hat{\mathcal{E}}^g_{0}(\tilde{\xi
}^n)\geq
\lim_{n}\hat{\mathcal{E}}^g_{0}(\zeta^n)= \hat{\mathcal
{E}}^g_{0}(\xi
)$, which is a contradiction.
Hence, $\liminf_n \hat{\mathcal{E}}_0^g ( \xi^n  )\geq
\hat
{\mathcal{E}}_0^g ( \xi )$.
\end{pf}
The preceding results allow us to derive a dual representation, by
means of the Fenchel--Moreau theorem, of the functional $\hat{\mathcal
{E}}^g(\cdot)$ at time zero.
%
%co4.10 #&#
\begin{cor}\label{cordualrepresentation}
Let $g$ be a generator fulfilling \textsc{(Pos)}, \textsc{(Lsc)} and
either \textsc{(Inc)} or \textsc{(Dec)}.
Assume that $g$ is jointly convex in $y$ and $z$.
Then, either $\hat{\mathcal{E}}^g_0\equiv+\infty$ or
%
%e4.31 #&#
%
\begin{equation}
\label{eqdual1} \hat{\mathcal{E}}^{g}_0(\xi)=
\mathcal{E}^{g}_0(\xi)=\sup_{\nu\in
L^\infty_+} \bigl
\{ E[\nu\xi]- \bigl(\hat{\mathcal{E}}^{g}_0
\bigr)^\ast ( \nu ) \bigr\},\qquad \xi\in L^1
\end{equation}
for the conjugate $(\hat{\mathcal{E}}^{g}_0)^\ast(\nu)=\sup_{\xi
\in
L^1} \{ E[\nu\xi]-\hat{\mathcal{E}}^{g}_0(\xi) \} $, where $\nu
\in
L^\infty$.
\end{cor}
\begin{pf}
Since $\hat{\mathcal{E}}^g_0>-\infty$ on $L^1$, either $\hat
{\mathcal
{E}}_0^g\equiv+\infty$ or $\hat{\mathcal{E}}_0^g$ is proper.
In the latter case, in view of Proposition~\ref{propproperties1} and
Theorem~\ref{thmmusc}, the function $\hat{\mathcal{E}}^{g}_0$ is
convex and $\sigma( L^1, L^\infty)$-lower semicontinuous on $L^1$.
Hence, the Fenchel--Moreau theorem yields the dual representation
\eqref{eqdual1}.
That the domain of $(\hat{\mathcal{E}}^{g}_0)^\ast$ is concentrated on
$L^\infty_+$ follows from the monotonicity of $\hat{\mathcal
{E}}^{g}_0$; see Proposition~\ref{propproperties1}.
\end{pf}
%
%re4.11 #&#
\begin{rem}\label{remlinkconvexriskmeasure}
Notice that if the generator in Corollary~\ref{cordualrepresentation}
does not depend on $y$, then by item (5) of
Proposition~\ref{propproperties1} the operator $\hat{\mathcal
{E}}^{g}_0(\cdot)$ is translation invariant.
Therefore, it is a lower semicontinuous, convex risk measure and
representation~\eqref{eqdual1} corresponds to the robust
representation of lower semicontinuous, convex risk measures; see F\"
{o}llmer and Schied \cite
{foellmer01}.
\end{rem}
Under additional integrability assumptions on the terminal condition we
may also formulate stability results for supersolutions in the set
$\mathcal{A}^1(\xi,g)$ introduced in~\eqref{setAp}.

%th4.12 #&#
\begin{thmm}\label{thmmmonconvfatouH1}
Suppose that the generator $g$ fulfills \textsc{(Pos)}, \textsc{(Lsc)},
\textsc{(Con)} and either \textsc{(Dec)} or \textsc{(Inc)}.
Let $(\xi^n)$ be a sequence in $L^0$, such that $\xi^n\geq\eta$, for
all $n\in\N$, where $(E[\eta| \mathcal{F}_{\cdot}])^*_{T}\in L^1$.
\begin{itemize}%[ font= \normalfont\itshape]
\item Suppose $(\xi^n)$ is increasing $P$-almost surely to
$\xi\in L^0$ and $\mathcal{A}^1(\xi,g)\neq\varnothing$. Then
$\mathcal
{E}^g_{t}(\xi)=\lim_{n}\mathcal{E}^g_{t}(\xi^n)$, $P$-almost surely,
for all $t\in[0,T]$.
\item Suppose $\mathcal{A}^1(\liminf_{n} \xi^n,g)\neq
\varnothing$. Then $\mathcal{E}^g_{t}(\liminf_{n} \xi^n)\leq\liminf_{n}\mathcal{E}^g_{t}(\xi^n)$, $P$-almost surely, for all $t\in[0,T]$.
\end{itemize}
\end{thmm}
We omit the proof of the preceding theorem, as it is a simple
adaptation of the proofs of Theorems \ref{thmmexistenceH1} and \ref
{thmmmonconvfatou}.
Note that Theorem~\ref{thmmmonconvfatouH1} is a weaker version of
Theorem~\ref{thmmmonconvfatou}.
Indeed, here, given a sequence $ ( \xi^n  )$ increasing to
$\xi$, we need to assume that $\mathcal{A}^1 (\xi,g )$ is
not empty.
The underlying reason being the lack of knowledge whether the limit
process $Y$, defined in the proof of Theorem~\ref{thmmmonconvfatou},
fulfills $Y^*_{T}\in L^1$.

The theorem above allows us to state the following result on $\mathcal
{L}^1$-lower semicontinuity of $\hat{\mathcal{E}}^g$.
Its proof is virtually the same as the proof of Theorem~\ref{thmmusc}.

%th4.13 #&#
\begin{thmm}\label{thmmlscH1}
Suppose that the generator $g$ fulfills \textsc{(Pos)}, \textsc{(Lsc)},
\textsc{(Con)} and either \textsc{(Dec)} or \textsc{(Inc)}.
Then $\xi\mapsto\hat{\mathcal{E}}^g_0 ( \xi )$ is
$\mathcal
{L}^1$-lower semicontinuous on its domain, that is, on
%
%e4.32 #&#
%
\begin{equation}
\label{setH1stability} \bigl\{ \xi\in L^0\dvtx \bigl(E\bigl[\xi^-|
\mathcal{F}_{\cdot}\bigr]\bigr)^*_{T} \in L^1\mbox
{ and }\mathcal{A}^1 ( \xi,g )\neq\varnothing \bigr\}.
\end{equation}
\end{thmm}
We conclude this section with a theorem on monotone stability with
respect to the generator.
%and $g$ satisfy \ref{cond02}.}
%
%th4.14 #&#
\begin{thmm}\label{thmmuscg}
Let $\xi\in L^0$ be a terminal condition, such that $\xi^-\in L^1$,
and let $(g^{n})$ be an increasing sequence of generators, which
converge pointwise to a generator $g$.
Suppose that each generator fulfills \textsc{(Pos)}, \textsc{(Lsc)},
\textsc{(Con)} and either \textsc{(Inc)} or \textsc{(Dec)}.
Then $\lim_{n}\hat{\mathcal{E}}^{g^n}_{0}(\xi)= \hat{\mathcal
{E}}^{g}_{0}(\xi)$.
If, in addition, $\lim_{n}\hat{\mathcal{E}}^{g^n}_{0}(\xi)<\infty$,
then $\mathcal{A}(\xi,g)\neq\varnothing$ and $\mathcal
{E}^{g^n}_{t}(\xi
)$ converges $P$-almost surely to $\mathcal{E}^{g}_{t}(\xi)$, for all
$t \in [ 0,T  ]$.
\end{thmm}
%
%re4.15 #&#
\begin{rem}\label{remgenstabext}
Under additional assumptions on the generators, one can prove
Fatou-type stability results for a $P\otimes dt$-almost sure converging
sequence of generators; see Gerdes, Heyne and Kupper \cite{GHK12} for details.
\end{rem}
\begin{pf*}{Proof of Theorem \ref{thmmuscg}}
Note that from Proposition~\ref{propproperties1}, we have $\hat
{\mathcal{E}}^{g^n}(\xi)\leq\hat{\mathcal{E}}^{g^{n+1}}(\xi)\leq
\cdots\leq\hat{\mathcal{E}}^{g}(\xi)$.
Hence, we may set $\hat Y_0=\lim_{n}\hat{\mathcal{E}}^{g^n}_0(\xi)$.
If $\hat Y_{0}=\infty$, then also $\hat{\mathcal{E}}^{g}_{0}(\xi
)=\infty$ and we are done.
Suppose that $\hat Y_{0}<\infty$.
By the same arguments as in the proof of Theorem~\ref
{thmmmonconvfatou}, we construct a c\`adl\`ag supermartingale $Y$.
With the same procedure as in Theorem~\ref{thmmmonconvfatou}, we
construct the candidate $\hat Z$.
It remains to show $(Y, \hat Z)\in\mathcal{A}(\xi,g)$.
However, this can be done similarly as in the proof of Theorem~\ref
{thmmexistence}.
We only show how to obtain the analogue of \eqref{thmmexistencestep54}.
Note first that the pointwise convergence of the generators implies
that $(g^k(Y,\hat Z))$ converges $P\otimes dt$-almost surely to
$g(Y,\hat Z)$.
Hence, Fatou's lemma yields
%
%e4.33 #&#
%e4.34 #&#
%
\begin{eqnarray}
\label{eqproofthmmuscgfirstfatou}&& Y_{s}-\int_{s}^{t}g_{u}(Y_u,
\hat Z_u)\,du+\int_{s}^{t}\hat
Z_u \,dW_u
\nonumber
\\[-8pt]
\\[-8pt]
\nonumber
&&\qquad\geq\limsup_{k} \biggl(Y_{s}-\int
_{s}^{t}g^k_{u}(Y_u,
\hat Z_u)\,du+\int_{s}^{t}\hat
Z_u \,dW_u \biggr).
\end{eqnarray}
As in the previous proof, we use the expression in the bracket on the
right-hand side to obtain
\begin{eqnarray*}
&& Y_{s}-\int_{s}^{t}g^k_{u}(Y_u,
\hat Z_u)\,du+\int_{s}^{t}\hat
Z_u \,dW_u
\\
&&\qquad\geq\limsup_{n}\sum_{i=n}^{M^{(n)}}
\lambda^{(n)}_{i} \biggl( Y_{s}^{i}-
\int_{s}^{t}g^k_{u}\bigl(
Y^{i}_u, Z_u^{i}\bigr)\,du+\int
_{s}^{t} Z_u^{i}\,dW_u
\biggr).
\end{eqnarray*}
Since on the right-hand side we consider the $\limsup$ with respect to
$n$ and $k$ being fixed for the moment, we may assume $k\leq n$, which
entails by monotonicity of the sequence of generators
\begin{eqnarray*}
&& Y_{s}-\int_{s}^{t}g^k_{u}(Y_u,
\hat Z_u)\,du+\int_{s}^{t}\hat
Z_u \,dW_u
\\
&&\qquad\geq\limsup_{n}\sum_{i=n}^{M^{(n)}}
\lambda^{(n)}_{i} \biggl( Y_{s}^{i}-
\int_{s}^{t}g^i_{u}\bigl(
Y^{i}_u, Z_u^{i}\bigr)\,du+\int
_{s}^{t} Z_u^{i}\,dW_u
\biggr).
\end{eqnarray*}
From here, we obtain as before $Y_{s}-\int_{s}^{t}g^k_{u}(Y_u,\hat
Z_u)\,du+\int_{s}^{t}\hat Z_u \,dW_u \geq Y_{t}$, where the right-hand side
does not depend on $k$ anymore.
Combined with \eqref{eqproofthmmuscgfirstfatou}, this yields the
analogue of \eqref{thmmexistencestep54}.
\end{pf*}

%
%s4.3 #&#
\subsection{Nonpositive generators}\label{secgnonpostive}
In this section we extend our results to generators that are not
necessarily positive.
%This is important with regards to applications in mathematical
%finance, where the generators are quite often of the linear-quadratic
%type.
%It turns out that we can extend the scope of our theorems to cover
%precisely some of these situations.
Using some measure change, the positivity assumption on the generator
$g$ can be relaxed to a linear bound below.
This leads to optimal solutions under $P$, where the admissibility is
required with respect to the related equivalent probability measure.
More precisely, we say in the following that a generator $g$ is

\begin{longlist}
\item[\textsc{(Lb)}] linearly bounded from below if there exist
adapted measurable $\R^{1\times d}$ and $\R$-valued processes $a$ and
$b$, respectively, such that $g(y,z)\geq az^\top+b$, for all $y,z\in
\R
\times\R^{1\times d}$.
Furthermore, $\int_{0}^tb_{s}\,ds\in L^1(P^{a})$, for all $t\in[0,T]$ and
\[
\frac{dP^{a}}{dP}=\mathcal{E} \biggl(\int a \,dW \biggr)_{T},
\]
defines an equivalent probability measure $P^{a}$.
\end{longlist}

%ex4.16 #&#
\begin{exep}
For instance, given a generator $g$, assume that there exists a
generator $\hat{g}$ independent of $y$ fulfilling \textsc{(Con)} and such
that $g \geq\hat{g}$.
Then, there exists an $\R^{1\times d}$-valued adapted measurable
process $a$ such that $g ( y,z  )\geq az^\top-\hat{g}^\ast
(a)$, for all $y,z\in\R\times\R^{1\times d}$, where $\hat{g}^\ast$
denotes the convex conjugate of $\hat g$.
%If $a$ defines a measure change and $b$ is integrable as above, then
%Assumption~\ref{cond00bis} holds.
\end{exep}
In the following, we say that $Z$ is $a$-admissible, if $\int_{}^{}Z
\,dW^a$ is a $P^{a}$-super\-martingale, where $W^a=(W^1-\int_{}^{}a^1
\,ds,\ldots, W^d-\int_{}^{}a^d \,ds )^\top$ is the respective Brownian
motion under $P^a$.
We are interested in the sets
%
%e4.35 #&#
%
\begin{equation}
\label{setAa}\quad  \mathcal{A}^{a} (\xi,g )= \bigl\{ (Y,Z) \in\mathcal {S}
\times \mathcal{L} \dvtx Z \mbox{ is $a$-admissible and } \eqref {eqcentralineq}
\mbox{ holds} \bigr\},
\end{equation}
and define the random process
%
%e4.36 #&#
%
\begin{equation}
\label{eqdefiphia}\quad  \hat{\mathcal{E}}^{g,a}_{t}( \xi)= \essinf
\bigl\{ Y_t\in L^0 (\mathcal {F}_t ) \dvtx (Y,Z )
\in\mathcal{A}^{a} (\xi,g ) \bigr\},\qquad t \in [ 0,T ].
\end{equation}
The analogue of Theorem~\ref{thmmexistence} is given as follows:
%
%th4.17 #&#
\begin{thmm}\label{thmmexistencegeng}
Let $g$ be a generator fulfilling \textsc{(Lb)}, \textsc{(Lsc)},
\textsc{(Con)} and either \textsc{(Inc)} or \textsc{(Dec)} and $\xi\in L^0$ be a
terminal condition, such that $\xi^-\in L^1(P^{a})$.
If $\mathcal{A}^{a}(\xi,g) \neq\varnothing$, then there exists a unique
minimal supersolution $(\hat Y,\hat Z)\in\mathcal{A}^{a}(\xi,g)$.
Moreover, $\mathcal{E}^g ( \xi )$ is the value process of
the minimal supersolution, that is, $(\mathcal{E}^g ( \xi
),\hat{Z})\in\mathcal{A}^{a}(\xi,g)$.
\end{thmm}
The analogues of Theorems \ref{thmmmonconvfatou} and~\ref
{thmmusc} read as follows.
%
%th4.18 #&#
\begin{thmm}\label{thmmmonconvfatougeng}\label{thmmuscgeng}
Suppose that the generator $g$ fulfills \textsc{(Lb)}, \textsc{(Lsc)}, \textsc{(Con)} and either \textsc{(Inc)} or \textsc{(Dec)}.
Let $(\xi^n)$ be a sequence in $L^0$, such that $\xi^n\geq\eta$, for
all $n\in\N$, where $\eta\in L^1(P^{a})$.
\begin{itemize}%[ font= \normalfont\itshape]
\item Monotone convergence: If $(\xi^n)$ is increasing
$P$-almost surely to $\xi\in L^0$, then $\hat{\mathcal
{E}}^{g,a}_{0}(\xi
)=\lim_{n}\hat{\mathcal{E}}^{g,a}_{0}(\xi^n)$.
\item Fatou's lemma: $\hat{\mathcal
{E}}^{g,a}_{0}(\liminf_{n} \xi^n)\leq\liminf_{n}\hat{\mathcal{E}}^{g,a}_{0}(\xi^n)$.
\end{itemize}
In particular, $\hat{\mathcal{E}}_0^{g,a}$ is $L^1(P^a)$-lower semicontinuous.
\end{thmm}
%
%As for the analogues of Theorem~\ref{thmmusc}.
% Given a generator $g$ fulfilling \ref{cond00bis}, \ref{cond01} and
%either \ref{cond02} or \ref{cond03}, then $\xi\mapsto
We only prove the first theorem.
\begin{pf*}{Proof of Theorem~\ref{thmmexistencegeng}}
In the setting of Section~\ref{secex}, given a positive generator
$\bar{g}$ and a random variable $\zeta$, let us denote by $\mathcal
{A} ( \zeta, \bar{g},W^a )$ the set defined in~\eqref{setA}
to indicate the dependence of this set on the Brownian motion $W^a$ and
the respective probability measure $P^a$.
Let us now define the generator $\bar g$ as
%
%e4.37 #&#
%
\begin{equation}
\label{eqdefnbarg}\qquad \bar{g} ( y,z )=g \biggl( y+\int_{0}^\cdot
b_{s}\,ds,z \biggr)-az^\top-b \qquad\mbox{for all } (y,z) \in\R
\times\R^{1\times d}.
\end{equation}
By assumption~\textsc{(Lb)}, this generator fulfills \textsc{(Pos)},
\textsc{(Lsc)}, \textsc{(Con)} and either \textsc{(Inc)} or \textsc{(Dec)}.
Since $\int_{}^{}Z \,dW^a$ is a $P^a$-supermartingale, a simple
inspection shows that the affine transformation $\bar{Y}=Y-\int_{}^{}b
\,ds$ and $\bar{Z}=Z$ yields a one-to-one relation between $\mathcal
{A}^a( \xi,g )$ and $\mathcal{A}( \xi-\int_{0}^{T}b_s \,ds,\bar{g},W^a)$.
Hence, the assumptions of Theorem~\ref{thmmexistence} are fulfilled
for $\bar{g}$ and $\mathcal{A}( \xi-\int_{0}^{T}b_s \,ds,\bar{g},W^a )$,
and thus its application completes the proof.
\end{pf*}

%re4.19 #&#
\begin{rem}\label{remH1nonpos}
Note that if $(E^{a}[ (\xi-\int_{0}^Tb_{s}\,ds)^- | \mathcal{F}_{\cdot
} ])^\ast_{T} \in L^1 ( P^a  )$, then Theorem~\ref
{thmmexistenceH1} applies in the same way; that is, under the
assumptions of Theorem~\ref{thmmexistencegeng}, if
\[
\mathcal{A}^{1,a} ( \xi,g ):=\bigl\{ ( Y,Z )\in \mathcal{A}^a
( \xi,g )\dvtx Z \in\mathcal{L}^1 \bigl( P^a \bigr) \bigr\}
\neq\varnothing,
\]
then $\mathcal{E}^{g,a} ( \xi )$ is the value process of the
minimal supersolution with unique control process $Z \in\mathcal
{L}^1 ( P^a  )$.
\end{rem}

\section*{Acknowledgments}
We thank Patrick Cheridito, Hans F\"ollmer, Ramon Van Handel, Ulrich
Horst and Reinhard
Schmidt for helpful comments and fruitful discussions.
We thank an anonymous referee for careful reading and helpful comments.

\printaddresses

\end{document}